\pdfoutput=1
\documentclass[11pt]{article}
\usepackage[margin=1.15in]{geometry}
\usepackage{amsmath,amssymb,amsthm}
\usepackage{tikz}
\usepackage{enumitem}
\usepackage[colorlinks=true,linkcolor=blue,citecolor=blue,urlcolor=blue]{hyperref}

\newtheorem{theorem}{Theorem}[section]
\newtheorem{lemma}[theorem]{Lemma}
\newtheorem{proposition}[theorem]{Proposition}
\newtheorem{corollary}[theorem]{Corollary}
\theoremstyle{definition}

\newtheorem{question}[theorem]{Question}
\theoremstyle{remark}
\newtheorem{remark}[theorem]{Remark}

\newcommand{\G}{\Gamma}
\newcommand{\g}{\gamma}
\newcommand{\pw}{\operatorname{pow}}
\newcommand{\dist}{\operatorname{dist}}
\newcommand{\Area}{\operatorname{Area}}

\begin{document}

\title{Tangent lengths and invariants of Poncelet polygons between two circles}

\author{Quang Hung Tran\\[2pt]
\small High School for Gifted Students, Vietnam National University, Hanoi, Vietnam\\
\small\texttt{tranquanghung@hus.edu.vn}}
\date{November 15, 2026}

\maketitle

\begin{abstract}
A polygon inscribed in one circle and circumscribed about another moves in a one-parameter family, by Poncelet's porism. We study this family through the tangent lengths from the vertices to the inner circle and find a number of quantities that remain constant during the motion. For example, the sum of the fourth powers of the distances from the circumcenter to the points of contact, plus twice the squared area, is constant. We also obtain a version of Euler's formula for $n$-gons, extend a problem on bicentric hexagons from \emph{Kvant}, and find a second quadratic invariant that goes along with the classical theorem of Weill. The basic invariants have elementary proofs; for the main ones we give two proofs, one of them without derivatives. Reciprocation at a limiting point of the two circles turns these results into statements about billiards in ellipses, and gives short proofs of theorems of Akopyan, Schwartz and Tabachnikov and of Bialy and Tabachnikov. Using the elliptic curve of the family, we show that the area determines the polygon up to a reflection, find the polygons of largest and smallest area, which proves two conjectures of Cie\'slak and Mozgawa, and find all points for which the product of the distances to the vertices is constant.
\end{abstract}



\section{Introduction}

\subsection{Poncelet's theorem}
Let $C$ and $D$ be two conics in the plane. Start at a point $A_1$ of $C$, draw a tangent line to $D$, let it meet $C$ again at $A_2$, draw the other tangent from $A_2$ to $D$, and so on. We assume that the tangents needed here are real, as is the case when $D$ lies inside $C$. Poncelet's porism \cite{Poncelet} says that if this process closes up after $n$ steps for one starting point, then it closes up after $n$ steps for every starting point. The closed polygons $A_1A_2\dots A_n$ form a one-parameter family, and we call each of them a \emph{Poncelet polygon}. For the history and for proofs we refer to \cite{BKOR,DelCentina,DR,Flatto,GH}.

Poncelet proved more. Recall that the \emph{pencil} spanned by $C$ and $D$ consists of the conics $\alpha C+\beta D=0$.

\begin{theorem}[Poncelet's general theorem]\label{thm:PGT}
Let $A_1A_2\dots A_n$ run over the family of Poncelet polygons inscribed in $C$ and circumscribed about $D$. Then for every $m$ with $1\le m\le n-1$ the diagonals $A_iA_{i+m}$ remain tangent to a fixed conic $D_m$ of the pencil spanned by $C$ and $D$. The conic $D_m$ may be degenerate, and then the diagonals pass through its singular point, which is what tangency means for such a conic, see Section \ref{ss:diag}.
\end{theorem}

In this paper both conics are circles, one inside the other. We write $\G$ for the outer circle, with center $O$ and radius $R$, and $\g$ for the inner one, with center $I$ and radius $r$, and we put $d=|OI|$, so that $d+r<R$. If $d=0$ the family consists of rotations of a regular polygon and there is nothing to say, so we assume $d>0$ throughout. A polygon inscribed in $\G$ and circumscribed about $\g$ is called \emph{bicentric}. For triangles the two circles are related by the formula of Chapple and Euler,
\begin{equation}\label{eq:euler}
d^2=R^2-2Rr ,
\end{equation}
and for quadrilaterals by the formula of Fuss \cite{Chapple,Fuss}. Jacobi \cite{Jacobi} wrote the vertices of a bicentric $n$-gon in terms of his elliptic functions, and this gives a closure condition for every $n$; see \cite{BKOR}.

\subsection{Pencils of circles and limiting points}
A pencil of circles is the family of circles that share a radical axis with two given circles. All centers lie on one line, perpendicular to the radical axis. There are three kinds of pencils \cite{AZ,CG}. If the two circles meet at two points, all circles of the pencil pass through these points (an \emph{elliptic} pencil). If they touch, all circles touch at the same point (a \emph{parabolic} pencil). If they are disjoint, no two circles of the pencil meet (a \emph{hyperbolic} pencil). Our circles $\G$ and $\g$ are disjoint, so we are in the third case.

A hyperbolic pencil contains exactly two circles of radius zero. These two points $L$ and $L'$ are the \emph{limiting points} of the pencil. They lie on the line of centers and are symmetric in the radical axis. We always denote by $L$ the limiting point that lies inside $\g$ and by $L'$ the one outside $\G$; the points $O$, $I$, $L$ lie on the line of centers in this order. If $x_0=(R^2+d^2-r^2)/(2d)$ is the distance from $O$ to the radical axis, then
\[
|OL|=x_0-\sqrt{x_0^2-R^2},\ \ \  |OL'|=x_0+\sqrt{x_0^2-R^2}.
\]
We shall use three standard facts.
\begin{enumerate}[label=(\alph*)]
\item $L$ and $L'$ are inverse to each other in every circle of positive radius of the pencil. Hence the polar of $L$ with respect to every such circle is the same line, namely the perpendicular to the line of centers at $L'$.
\item An inversion centered at $L$ sends the circles of the pencil to concentric circles.
\item By Theorem \ref{thm:PGT}, for $1\le m<n/2$ the diagonals $A_iA_{i+m}$ of a bicentric polygon touch a circle $\g_m$ of the pencil of $\G$ and $\g$, see Section \ref{ss:diag}. When $n=2N$ is even and $m=N$, the conic of Theorem \ref{thm:PGT} degenerates to the point $L$, so the main diagonals $A_iA_{i+N}$ all pass through $L$; see \cite{GOR,KM} and Proposition \ref{prop:conj}.
\end{enumerate}

\subsection{From two circles to confocal conics}\label{ss:intro-polar}
Take a circle $\omega$ with center $P$ and radius $\rho$. Reciprocation in $\omega$ sends a point $X\neq P$ to its polar line, the line perpendicular to $PX$ at distance $\rho^2/|PX|$ from $P$, and sends lines to points accordingly. The reciprocal of a circle with center $Z$ and radius $\varrho$ is a conic with focus $P$ and eccentricity $|PZ|/\varrho$ \cite{AZ}.

Now let $P$ be a limiting point of our pencil. Then the reciprocals of $\G$ and $\g$ are not just two conics with a common focus, they are \emph{confocal}: the second focus is common as well \cite{RGR}. If $P=L$, both eccentricities are less than $1$ and we get two confocal ellipses, $E=\g^{\circ}$ outside and $E'=\G^{\circ}$ inside. If $P=L'$, both eccentricities exceed $1$ and we get two confocal hyperbolas. A bicentric polygon turns into a polygon inscribed in $E$ whose sides touch $E'$, and such a polygon is a periodic billiard trajectory in $E$. So the Poncelet family between two circles and the family of periodic billiard trajectories in an ellipse are two pictures of one object. This was the point of view of Roitman, Garcia and Reznik in \cite{RGR}.

The correspondence preserves incidences but not lengths, so a metric invariant on one side does not automatically give a metric invariant on the other. Each side has its own list. On the billiard side the list is long: starting with computer experiments of Reznik \cite{RGK1,RGK2}, many conserved quantities were found, and several were proved by Akopyan, Schwartz and Tabachnikov \cite{AST}, Bialy and Tabachnikov \cite{BT}, Chavez-Caliz \cite{Chavez} and others. On the side of two circles the classical result is the theorem of Weill \cite{Weill1878}: the centroid of the points of contact is a fixed point. Recently \cite{RGR} added the sum of the cosines of the angles of the polygon, and \cite{KM} found invariants of bicentric hexagons.

\subsection{What is done here}
The present paper started from the following problem, which the author proposed in a short note in April 2026.

\medskip
\noindent\emph{Let $B_i$ be the point where the side $A_iA_{i+1}$ of a bicentric polygon touches $\g$, and let $S$ be the area of the polygon. Then $|OB_1|^4+\dots+|OB_n|^4+2S^2$ does not depend on the polygon.}
\medskip

Looking for a proof we were led to describe the polygon by the tangent lengths
\[
t_i=|A_iB_{i-1}|=|A_iB_i| ,
\]
as Radi\'c did in \cite{Radic05,Radic14}. In these variables everything becomes simple: two consecutive tangent lengths satisfy a symmetric quadratic relation, the distance from a vertex to a limiting point is proportional to $t_i$, and Weill's theorem says that $\sum t_it_{i+1}$ is constant. A second fact of the same kind, that $\sum_{i<j}t_it_j$ is constant, then gives the problem above and a good deal more.

Section \ref{s:results} states the results. Section \ref{s:tools} collects the tools, all of them classical, and Section \ref{s:proofs} contains the proofs. They use derivatives in two places, in Lemma \ref{lem:kin} and in the identity \eqref{eq:dW}. Section \ref{s:nocalc} shows how to avoid the first of them, which is the one needed for the main invariants. In Section \ref{s:conics} we translate the results to confocal ellipses and hyperbolas and say which translations we find attractive and which not. Section \ref{s:found} explains how the statements were found and gives further results which we can prove only with an elliptic curve: more invariants, the fact that the area determines the polygon up to symmetry, the polygons of extreme area and perimeter, and the list of all points for which the product of the distances to the vertices is constant. It also contains a theorem on the tangents of the half-angles at $L$ and a recursion for the constants $\lambda_m$.

\section{Results}\label{s:results}

\subsection{Notation}
Throughout, $\mathcal P=A_1A_2\dots A_n$ is a convex polygon inscribed in $\G$ and circumscribed about $\g$, labelled counterclockwise, with indices taken modulo $n$. The side $A_iA_{i+1}$ touches $\g$ at $B_i$, the tangent lengths are $t_i=|A_iB_i|=|A_iB_{i-1}|$, and $S$ is the area of $\mathcal P$. Since $\g$ is the incircle,
\begin{equation}\label{eq:S}
|A_iA_{i+1}|=t_i+t_{i+1},\ \ \  S=r\,(t_1+\dots+t_n).
\end{equation}
We use the constants
\[
k=R^2-d^2,\ \ \  c=R^2-d^2-r^2,\ \ \  q=\sqrt{\bigl((R-d)^2-r^2\bigr)\bigl((R+d)^2-r^2\bigr)} ,
\]
so that $q^2=c^2-4d^2r^2$. Note that $q$ is the product of the smallest and the largest possible tangent length. We also put
\[
W=\sum_{i=1}^n t_it_{i+1},\ \ \  W_m=\sum_{i=1}^n t_it_{i+m},\ \ \  e_2=\sum_{i<j}t_it_j .
\]
For the limiting points we write $\mu=\sqrt{|OL|/|OI|}$ and $\mu'=\sqrt{|OL'|/|OI|}$. Finally, for $1\le m<n/2$, $\g_m$ is the circle of Theorem \ref{thm:PGT} touched by the diagonals $A_iA_{i+m}$, with center $I_m$, and $\lambda_m=\sqrt{|OI_m|/|OI|}$; thus $\g_1=\g$ and $\lambda_1=1$.

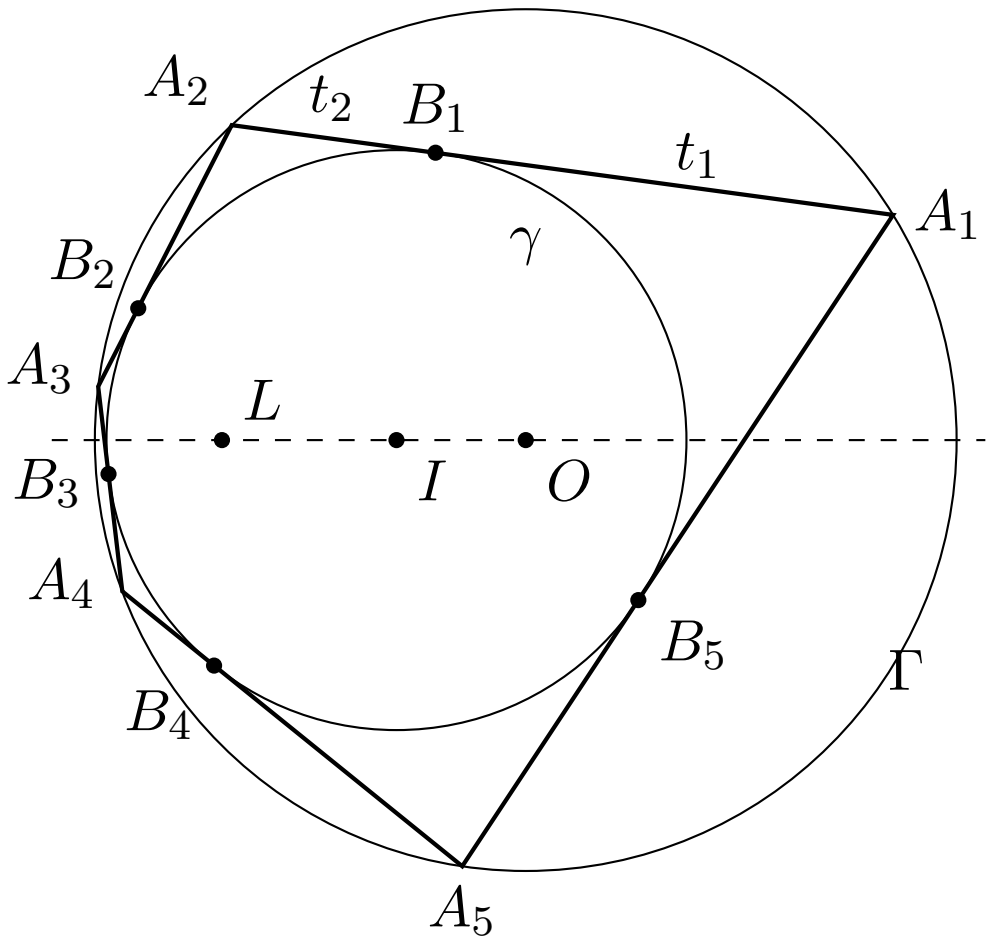
\begin{figure}[ht]
\centering
\begin{tikzpicture}[scale=0.95]
\draw (0.9,0) circle (3.0);
\draw (0,0) circle (2.019);
\draw[thick] (3.458,1.568)--(-1.147,2.193)--(-2.077,0.373)--(-1.909,-1.053)--(0.458,-2.967)--cycle;
\foreach \p in {(0.272,2.001),(-1.798,0.918),(-2.005,-0.236),(-1.270,-1.570),(1.684,-1.114)} \fill \p circle (1.6pt);
\fill (0.9,0) circle (1.6pt) node[below right] {$O$};
\fill (0,0) circle (1.6pt) node[below right] {$I$};
\fill (-1.215,0) circle (1.6pt) node[above right] {$L$};
\draw[dashed] (-2.4,0)--(4.1,0);
\node[right] at (3.458,1.568) {$A_1$};
\node[above left] at (-1.147,2.193) {$A_2$};
\node[left] at (-2.1,0.5) {$A_3$};
\node[left] at (-1.95,-1.0) {$A_4$};
\node[below] at (0.458,-2.967) {$A_5$};
\node[above] at (0.272,2.001) {$B_1$};
\node[above left] at (-1.798,0.918) {$B_2$};
\node[left] at (-2.05,-0.3) {$B_3$};
\node[below left] at (-1.270,-1.60) {$B_4$};
\node[below right] at (1.684,-1.114) {$B_5$};
\node at (2.1,1.98) {$t_1$};
\node at (-0.45,2.38) {$t_2$};
\node at (3.55,-1.6) {$\G$};
\node at (0.9,1.35) {$\g$};
\end{tikzpicture}
\caption{A bicentric pentagon, its points of contact and the limiting point $L$.}\label{fig:one}
\end{figure}

\subsection{Identities of the first degree}
The first group of results holds for each bicentric polygon separately. Nothing moves here, but the right-hand sides depend only on the two circles, so these are invariants of the Poncelet family as well.

\begin{theorem}\label{thm:first}
For every bicentric $n$-gon the following hold.
\begin{enumerate}[label=\textup{(\alph*)}]
\item The short diagonals satisfy
\begin{equation}\label{eq:diag2}
\sum_{i=1}^n |A_iA_{i+2}|=\frac{4R}{R^2-d^2}\,S .
\end{equation}
Equivalently, if $\Delta$ denotes the ratio of $\sum|A_iA_{i+2}|$ to the perimeter $\sum|A_iA_{i+1}|$, then
\begin{equation}\label{eq:euler-n}
\frac{1}{R-d}+\frac{1}{R+d}=\frac{\Delta}{r},\ \ \text{that is,}\ \  d^2=R^2-\frac{2Rr}{\Delta}.
\end{equation}
\item The polygon of contact points satisfies $\Area(B_1B_2\dots B_n)=\dfrac{r^2}{R^2-d^2}\,S$.
\item For every $m$ with $1\le m<n/2$,
\begin{equation}\label{eq:diagm}
\sum_{i=1}^n|A_iA_{i+m}|=\frac{2\lambda_m}{r}\,S .
\end{equation}
If $n=2N$ is even, then for the main diagonals
\begin{equation}\label{eq:kvant}
\frac{S}{|A_1A_{N+1}|+|A_2A_{N+2}|+\dots+|A_NA_{2N}|}=\frac{r}{\mu}=r\sqrt{\frac{|OI|}{|OL|}} .
\end{equation}
\item $\displaystyle\sum_{i=1}^n|A_iA_{i+1}|\cdot|OB_i|^2=(d^2+r^2)\sum_{i=1}^n|A_iA_{i+1}|$.
\end{enumerate}
\end{theorem}

For a triangle $A_iA_{i+2}$ is a side, so $\Delta=1$ and \eqref{eq:euler-n} is Euler's formula \eqref{eq:euler}. This is why we regard \eqref{eq:euler-n} as a form of Euler's formula for $n$-gons. It is not a closure condition, because $\Delta$ is determined by $R$, $r$, $d$; rather it says that $\Delta$ measures how far the pair of circles is from the triangle case. The proof of (a) gives more, namely $|A_iA_{i+2}|=\frac{2Rr}{R^2-d^2}\,(t_i+t_{i+2})$ for each $i$. Hence $\lambda_2=2Rr/(R^2-d^2)$ for $n\ge5$, and the center of the circle $\g_2$ touched by the short diagonals is given by
\[
|OI_2|=\frac{4R^2r^2d}{(R^2-d^2)^2}.
\]
We keep the notation $\lambda_2=2Rr/(R^2-d^2)$ for $n=3$ and $n=4$ as well. For $n=3$ this means $\g_2=\g$, and for $n=4$ the circle $\g_2$ is the point $L$ and $\lambda_2=\mu$.
Part (b) is due to Weill \cite{Weill1880}. For the hexagon, \eqref{eq:kvant} is Problem M2817 of \emph{Kvant} \cite{Kvant}; part (c) extends it to all $n$ and all $m$, and part (a) is the one case where the constant is an explicit rational function of $R$ and $d$ alone. For a quadrilateral $ABCD$ both \eqref{eq:diag2} and \eqref{eq:kvant} give $S/(|AC|+|BD|)=(R^2-d^2)/(2R)$.

\begin{remark}
For odd $n$ one may be tempted to replace the main diagonals by the segments joining each vertex to the point of contact on the opposite side. This does not work. In general these segments do not pass through $L$ (for a triangle they meet at the Gergonne point), and the ratio of $S$ to their sum is not constant: for $R=1$, $d=0.3$ and $n=5$ it varies from about $0.25616$ to $0.25637$.
\end{remark}

\subsection{Weill's theorem and its faces}
From now on the polygon moves.

\begin{theorem}[Weill and equivalent statements]\label{thm:weill}
The quantity $W=\sum t_it_{i+1}$ is constant over the Poncelet family. Moreover, the following relations hold for every bicentric polygon. Hence each of the quantities on their left-hand sides is constant as well, and the constancy of any one of them is equivalent to that of $W$ (in (iii) one takes a point $X$ with $\overrightarrow{IX}\cdot\overrightarrow{IO}\ne0$, for instance $X=O$).
\begin{enumerate}[label=\textup{(\roman*)}]
\item The centroid of $B_1,\dots,B_n$ lies on the line $OI$, at signed distance $(W-nc)/(2dn)$ from $I$, counted positive towards $O$. So it is a fixed point.
\item $\sum|OB_i|^2=nR^2-W$.
\item For every point $X$, the sum of the signed distances from $X$ to the lines of the sides, counted positive on the side of $I$, equals $nr-\dfrac{W-nc}{2d^2r}\,\overrightarrow{IX}\cdot\overrightarrow{IO}$. For $X=O$ this is $\bigl(n(R^2-d^2+r^2)-W\bigr)/(2r)$.
\item $\sum 1/|IA_i|^2=(n-W/k)/r^2$. Equivalently, the sum of the cosines of the interior angles of $\mathcal P$ equals $2W/k-n$.
\item $\sum|B_iB_{i+1}|^2=4r^2W/k$.
\item $\displaystyle\sum_{i=1}^n\frac{1}{|B_iA_i|\cdot|B_iA_{i+1}|}=\frac{W}{q^2}$ and $\displaystyle\sum_{i=1}^n\Bigl(\frac{t_i}{t_{i+1}}+\frac{t_{i+1}}{t_i}\Bigr)=\frac{2(nk-W)}{r^2}$.
\end{enumerate}
\end{theorem}

Statement (i) is the original theorem of Weill \cite{Weill1878}; see also \cite{Casey}. For a triangle (iii) with $X=O$ is Carnot's theorem, the sum being $R+r$. The invariance of the sum of cosines in (iv) was found and proved with elliptic functions in \cite{RGR}; the theorem above shows that it is Weill's theorem in another form. The second formula in (vi) is the one needed for the angles at the limiting points below.

\subsection{A second quadratic invariant}

\begin{theorem}\label{thm:second}
For every $m$ with $1\le m\le n-1$ the sum $W_m=\sum_i t_it_{i+m}$ is constant over the Poncelet family. Consequently $e_2=\sum_{i<j}t_it_j$ is constant, and so are
\[
S^2-r^2\sum_{i=1}^n|IA_i|^2,\ \ \  S^2-\frac{r^2}{2}\sum_{i=1}^n|A_iA_{i+1}|^2,\ \ \  S^2-2nr^2\,\overrightarrow{IO}\cdot\overrightarrow{IG},
\]
where $G$ is the centroid of the vertices. If $n=2N$, each product $t_it_{i+N}$ is separately constant and equal to $q$.
\end{theorem}

The last expression says that $S^2$ is an affine function of the projection of $G$ to the line of centers. By a theorem of Schwartz and Tabachnikov \cite{ST} the point $G$ moves on a circle, so this ties the area to the position of $G$ on that circle.

\subsection{Degree four and degree six}

\begin{theorem}\label{thm:fourth}
Over the Poncelet family,
\begin{align}
\sum_{i=1}^n|OB_i|^4+2S^2&=n(R^4-q^2)-2d^2W+4r^2e_2,\label{eq:four}\\
\sum_{i=1}^n|OB_i|^6+6\,d^2S^2&=\text{const}.\label{eq:six}
\end{align}
In terms of the products $p_i=t_it_{i+1}=|B_iA_i|\cdot|B_iA_{i+1}|$,
\[
\sum p_i^2+2S^2=4r^2e_2+2kW-nq^2,\ \ \  \sum p_i^3+6kS^2=12kr^2e_2+(4k^2-q^2+2r^4)W-2knq^2 .
\]
\end{theorem}

The coefficients $2$ and $6d^2$ on the left do not depend on $n$. Since $|OB_i|^2-d^2=k-p_i$, the theorem can also be read as follows: $\sum(|OB_i|^2-d^2)^2+2S^2$ and $\sum(|OB_i|^2-d^2)^3$ are constant. In degree eight this stops: $\sum|OB_i|^8$ is no longer a combination of $S^2$ and a constant. In degree ten there is again a simple invariant, see Proposition \ref{prop:ten} in Section \ref{s:found}.

\subsection{Products}

\begin{theorem}\label{thm:product}
Let $s$ be a smooth parameter on the Poncelet family. Then
\begin{equation}\label{eq:dW}
\frac{\partial}{\partial s}\sum_{i=1}^nt_it_{i+1}=(R^2-d^2)\,\frac{\partial}{\partial s}\log\prod_{i=1}^n|IA_i|^2 .
\end{equation} Hence $|IA_1|\cdot|IA_2|\cdots|IA_n|$ is constant. More generally $\prod_i|I_mA_i|$ is constant for $1\le m<n/2$, and $t_1t_2\cdots t_n=q^{n/2}$ if $n$ is even. If $n$ is odd and $m=(n-1)/2$, then
\begin{equation}\label{eq:odd-product}
\frac{\prod_{i=1}^n|A_iA_{i+1}|}{\prod_{i=1}^nt_i}=\frac{(2k_m)^n}{\prod_{i=1}^n|I_mA_i|^2},\ \ \ k_m=R^2-|OI_m|^2,
\end{equation}
so for odd $n$ the ratio of the product of the sides to the product of the tangent lengths is constant.
\end{theorem}

For a triangle these are the formulas $|IA|\cdot|IB|\cdot|IC|=4Rr^2$ and $abc/\bigl((s-a)(s-b)(s-c)\bigr)=4R/r$. All points $X$ for which $\prod|XA_i|$ is constant are listed in Theorem \ref{thm:dist}.

\subsection{The view from the limiting points}
Let $\varphi_i=\angle A_iLA_{i+1}$ be the angle under which the side $A_iA_{i+1}$ is seen from $L$ and $\psi_i=\angle B_iLB_{i+1}$ the angle under which two consecutive points of contact are seen. Let $\varphi_i'$, $\psi_i'$ be the same angles at $L'$.

\begin{theorem}\label{thm:limit}
\begin{enumerate}[label=\textup{(\alph*)}]
\item $|LA_i|=\mu\,t_i$ and $|L'A_i|=\mu'\,t_i$ for all $i$. The line $LB_i$ bisects the angle $A_iLA_{i+1}$, and $L'B_i$ bisects the angle $A_iL'A_{i+1}$.
\item The sum of cosines $\sum\cos\varphi_i$ is constant:
\[
\sum_{i=1}^n\cos\angle A_iLA_{i+1}=\frac{(\mu^2-1)(nk-W)}{\mu^2r^2}-\frac{n}{\mu^2},
\]
and the same holds at $L'$ with $\mu'$ in place of $\mu$.
\item For $1\le m<n/2$ the sums $\sum_i\cos\angle A_iLA_{i+m}$ and $\sum_i\cos\angle A_iL'A_{i+m}$ are constant.
\item Let $\mathcal P^*$ be the conjugate polygon of Proposition \ref{prop:conj} and $S^*$ its area. Then
\[
\sum_{i=1}^n\cos^2\varphi_i+\frac{(\mu^2-1)^2}{2\mu^4r^4}\,\bigl(S^2+S^{*2}\bigr)
\]
is constant. If $n$ is even, then $S^*=S$, so that $\sum\cos^2\varphi_i+\dfrac{(\mu^2-1)^2}{\mu^4r^4}\,S^2$ is constant.
\item $\sum\cos^2\psi_i$ is constant, and at $L'$ the sum $\sum\cos^2\psi_i'$ is constant for $n\neq4$. Moreover, for every bicentric polygon,
\[
\prod_{i=1}^n|\cos\psi_i|\cdot\prod_{i=1}^n|IA_i|^2=\Bigl(\frac{k^2(\mu^2-\lambda_2^2)}{\mu^2}\Bigr)^{n/2},
\]
where $\lambda_2=2Rr/k$. So $\prod|\cos\psi_i|$ is constant by Theorem \ref{thm:product}, and then so is $\prod\cos\psi_i$, because on a connected family a cosine that does not vanish keeps its sign. The same holds at $L'$.
\end{enumerate}
\end{theorem}

The sum $\sum\cos^2\varphi_i$ itself is never constant; by part (d) this means that $S^2+S^{*2}$ is not constant, which is proved at the end of Section \ref{ss:curve}. The variation can be small: in several families with odd $n$ that we tried it showed only in the fifth or sixth digit, which makes numerical experiments treacherous.

Table \ref{tab:summary} lists the main invariants.

\begin{table}[ht]
\centering
\renewcommand{\arraystretch}{1.35}
\begin{tabular}{lll}
\hline
Quantity & Status & Where \\
\hline
$\sum|A_iA_{i+2}|\,/\,S=4R/(R^2-d^2)$ & identity & Thm.\ \ref{thm:first}(a)\\
$\sum|A_iA_{i+m}|\,/\,S$ & identity & Thm.\ \ref{thm:first}(c)\\
$\sum|OB_i|^2$, centroid of the $B_i$ & Weill & Thm.\ \ref{thm:weill}\\
$S^2-r^2\sum|IA_i|^2$ & new & Thm.\ \ref{thm:second}\\
$\sum|OB_i|^4+2S^2$,\ \ $\sum|OB_i|^6+6d^2S^2$ & new & Thm.\ \ref{thm:fourth}\\
$\prod|IA_i|$ & cf.\ \cite{BT} & Thm.\ \ref{thm:product}\\
$\sum\cos\angle A_iLA_{i+1}$ & cf.\ \cite{AST} & Thm.\ \ref{thm:limit}(b)\\
$\sum\cos^2\angle A_iLA_{i+1}+\text{const}\cdot(S^2+S^{*2})$ & new & Thm.\ \ref{thm:limit}(d)\\
$\sum\tan\frac12\angle A_iLA_{i+1}\cdot\tan\frac12\angle A_{i+1}LA_{i+2}$,\ \ $n\ne4$ & new & Thm.\ \ref{thm:tan}\\
$\sum(|OB_i|^2-d^2)^5-10(R^2-d^2)r^4S^2$ & new & Prop.\ \ref{prop:ten}\\
$e_{2j}(t_1,\dots,t_n)$,\ \ $e_{2j+1}(t_1,\dots,t_n)/S$ & $n=6$: \cite{KM} & Thm.\ \ref{thm:poly}\\
$\prod|XA_i|$ at $O$, $I_m$, $I_m'$ (and $L$, $L'$ for even $n$) & complete list & Thm.\ \ref{thm:dist}\\
\hline
\end{tabular}
\caption{Invariants of the Poncelet family between two circles.}\label{tab:summary}
\end{table}

\section{Tools}\label{s:tools}

\subsection{Powers in a pencil}
The power of a point $X$ with respect to a circle with center $Z$ and radius $\varrho$ is $\pw(X)=|XZ|^2-\varrho^2$. For a point outside the circle it is the square of the tangent length. Take coordinates with origin $O$ and the first axis along $OI$. Every circle of the pencil of $\G$ and $\g$ has an equation
\begin{equation}\label{eq:pencil}
x^2+y^2-R^2-\nu\,(x-x_0)=0,
\end{equation}
where $x=x_0$ is the radical axis and the center is $(\nu/2,0)$. The left-hand side is the power of the point $(x,y)$. Hence:

\begin{lemma}[Casey]\label{lem:casey}
Let $\g'$ be a circle of the pencil of $\G$ and $\g$, with center $J$. For every point $X$ of $\G$,
\[
\pw_{\g'}(X)=\frac{\overline{OJ}}{\overline{OI}}\ \pw_{\g}(X),
\]
where the bars denote signed lengths on the line of centers. More generally, for every point $X$ of the plane,
\begin{equation}\label{eq:three-powers}
\pw_{\g'}(X)-\lambda^2\pw_\g(X)=(1-\lambda^2)\,\pw_\G(X),\ \ \  \lambda^2=\frac{\overline{OJ}}{\overline{OI}} .
\end{equation}
\end{lemma}

\begin{proof}
By \eqref{eq:pencil}, $\pw_\nu(X)=\pw_\G(X)-\nu(x-x_0)$ for the circle with parameter $\nu$. Eliminating $x-x_0$ between the equations for $\g$ and $\g'$ gives \eqref{eq:three-powers}, and for $X$ on $\G$ the right-hand side vanishes.
\end{proof}

So for a circle $\g'$ of the pencil with $\overline{OJ}/\overline{OI}>0$ the tangents from the points of $\G$ are real, and their lengths are $\lambda$ times the tangent lengths to $\g$. This is in Casey's book \cite{Casey}. The lemma applies to the point circles $L$ and $L'$ too, and gives the first part of Theorem \ref{thm:limit}(a):
\begin{equation}\label{eq:LA}
|LA_i|=\mu\,t_i,\ \ \  |L'A_i|=\mu'\,t_i .
\end{equation}

\subsection{Tangent lengths}
Take now the origin at $I$ and let $O=(d,0)$. Write $B_i=r(\cos\theta_i,\sin\theta_i)$.

\begin{lemma}\label{lem:tt}
For every $i$,
\begin{equation}\label{eq:tt}
t_i-t_{i+1}=2d\sin\theta_i,\ \ \  t_it_{i+1}=c+2dr\cos\theta_i,\ \ \  |OB_i|^2=R^2-t_it_{i+1}.
\end{equation}
\end{lemma}

\begin{proof}
Let $u=(\cos\theta_i,\sin\theta_i)$ and $v=(-\sin\theta_i,\cos\theta_i)$. The points of the side line are $ru+sv$, and $A_i$, $A_{i+1}$ correspond to $s=-t_i$ and $s=t_{i+1}$. The condition $|ru+sv-O|^2=R^2$ reads $s^2+2d\sin\theta_i\cdot s-(c+2dr\cos\theta_i)=0$. The sum and the product of the roots give the first two formulas. The third is the power of $B_i$ with respect to $\G$.
\end{proof}

Eliminating $\theta_i$ we get a relation between neighbours,
\begin{equation}\label{eq:Phi}
r^2(t_i-t_{i+1})^2+(t_it_{i+1}-c)^2=4d^2r^2 ,
\end{equation}
which we shall use in the form
\begin{equation}\label{eq:p2}
p_i^2=-r^2(t_i^2+t_{i+1}^2)+2k\,p_i-q^2,\ \ \  p_i=t_it_{i+1}.
\end{equation}
Relation \eqref{eq:Phi} is symmetric and of degree two in each variable. In a slightly different form it appears in the work of Radi\'c \cite{Radic05}.

\begin{lemma}\label{lem:vieta}
For every $i$,
\begin{equation}\label{eq:vieta}
t_{i-1}+t_{i+1}=\frac{2(R^2-d^2)\,t_i}{t_i^2+r^2}.
\end{equation}
\end{lemma}

\begin{proof}
Let $\alpha$ be the polar angle of $A_i$ seen from $I$ and $\phi=\angle A_iIB_i$, so that $\theta_{i-1}=\alpha-\phi$, $\theta_i=\alpha+\phi$, $\sin\phi=t_i/|IA_i|$, $\cos\phi=r/|IA_i|$ and $|IA_i|^2=t_i^2+r^2$. By Lemma \ref{lem:tt},
\[
t_{i-1}+t_{i+1}-2t_i=2d\bigl(\sin(\alpha-\phi)-\sin(\alpha+\phi)\bigr)=-4d\cos\alpha\sin\phi ,
\]
\[
t_i(t_{i-1}+t_{i+1})=2c+2dr\bigl(\cos(\alpha-\phi)+\cos(\alpha+\phi)\bigr)=2c+4dr\cos\alpha\cos\phi .
\]
Multiply the first equation by $r^2$, the second by $t_i$ and add. The terms with $\cos\alpha$ cancel and we are left with $(t_{i-1}+t_{i+1})(t_i^2+r^2)=2(c+r^2)t_i$.
\end{proof}

\subsection{Diagonals}\label{ss:diag}
We use Theorem \ref{thm:PGT} in the following form: for $1\le m<n/2$ the diagonals $A_iA_{i+m}$ touch a fixed circle $\g_m$ of the pencil, $\g_m$ lies inside $\G$, and hence the point of contact lies between $A_i$ and $A_{i+m}$. Theorem \ref{thm:PGT} only provides a conic of the pencil, so let us explain why it is a circle inside $\G$.

Let $f$ be the Poncelet map of $\G$, so that $f(A_i)=A_{i+1}$, and let $\sigma$ be the reflection in the line $OI$; then $\sigma f\sigma=f^{-1}$. The pencil of conics spanned by $\G$ and $\g$ consists of the circles \eqref{eq:pencil} and of three degenerate conics: the point circles $L$ and $L'$, and the pair formed by the radical axis and the line at infinity. Each of them has rank two. If $M$ is the matrix of such a conic, we call a line with the coordinate vector $\ell$ tangent to it when $\ell^{\mathsf T}\operatorname{adj}(M)\,\ell=0$, the equation of the dual conic; this means that the line passes through the singular point, which is $L$, $L'$, or the point at infinity of the radical axis. So if the conic $D_m$ of Theorem \ref{thm:PGT} were degenerate, all the lines $Xf^m(X)$, $X\in\G$, would pass through one of these three points. If the point is $L$ or $L'$, then $f^m$ is an involution, so $f^{2m}=\mathrm{id}$, which is impossible for $m<n/2$. If it is the point at infinity, then $f^m=\sigma$, while $f^m$ preserves the orientation of $\G$. So the conic is one of the circles \eqref{eq:pencil}. Take the origin at $O$ and let $(x,0)$ be its center and $\varrho$ its radius; by \eqref{eq:pencil}, $\varrho^2=x^2-2x_0x+R^2$, and $x_0>R$. Real lines touch the circle, so $\varrho^2>0$, that is, $x<|OL|$ or $x>|OL'|$. The map $\sigma f^m$ reverses the orientation of $\G$, so it has a fixed point $A$, and $f^m(A)=\sigma(A)\ne A$. Thus one of the diagonals is a chord of $\G$ perpendicular to $OI$. If $x<0$, the points of $\G$ have negative power with respect to the circle by Lemma \ref{lem:casey}, and there are no tangents through them. If $x>|OL'|$, then $x-\varrho>x_0>R$, and the tangents to the circle perpendicular to $OI$ do not meet $\G$. The value $x=0$ gives $\G$ itself. There remains $0<x<|OL|$, and then $\varrho<R-x$, so the circle lies inside $\G$.

Moreover, $\g_m$ lies on the left of every chord directed from $A_i$ to $A_{i+m}$. To see this, let $d_m=|OI_m|=x$ and let $r_m$ be the radius of $\g_m$. Since $|OL|$ is a root of $z^2-2x_0z+R^2=0$,
\[
r_m^2-|LI_m|^2=d_m^2-2x_0d_m+R^2-(|OL|-d_m)^2=2(x_0-|OL|)(|OL|-d_m)>0,
\]
so $L$ lies inside $\g_m$. For a point $X$ of $\G$ let $\alpha(X)\in(0,2\pi)$ be the angle of the counterclockwise rotation about $L$ that takes the ray $LX$ to the ray $Lf^m(X)$. The chord $Xf^m(X)$ touches $\g_m$, so it does not pass through $L$, and $\alpha(X)\ne\pi$. As $\alpha$ is continuous on $\G$, its values lie on one side of $\pi$. The polygon contains $L$, and the rays $LA_1,\dots,LA_n$ follow each other counterclockwise, so $\sum_i\alpha(A_i)=2m\pi<n\pi$. Hence $\alpha<\pi$ everywhere, that is, $L$, and with it $\g_m$, lies on the left of the chord directed from $A_i$ to $A_{i+m}$.

By Lemma \ref{lem:casey} the tangent length from $A_i$ to $\g_m$ is $\lambda_mt_i$, so
\begin{equation}\label{eq:diag-len}
|A_iA_{i+m}|=\lambda_m(t_i+t_{i+m}).
\end{equation}
For $n=2N$ and $m=N$ the diagonals pass through $L$ (see the Introduction, or Proposition \ref{prop:conj} below), which lies between $A_i$ and $A_{i+N}$, and \eqref{eq:diag-len} holds with $\lambda_N=\mu$ by \eqref{eq:LA}.

\subsection{The conjugate polygon}
\begin{proposition}\label{prop:conj}
Let $A_i^*$ be the second point of intersection of the line $A_iL$ with $\G$. Then $\mathcal P^*=A_1^*A_2^*\dots A_n^*$ is again a Poncelet polygon for $\G$ and $\g$, and its tangent lengths are
\[
t_i^*=\frac{q}{t_i}.
\]
If $n=2N$ then $A_i^*=A_{i+N}$.
\end{proposition}

\begin{proof}
Let $\ell'$ be the common polar of $L$ with respect to the circles of the pencil, see fact (a) of the Introduction, and let $h$ be the harmonic homology with center $L$ and axis $\ell'$: it fixes $L$ and the points of $\ell'$, and sends $X$ to the point $X^h$ of the line $LX$ such that $X$, $X^h$ are harmonic conjugates with respect to $L$ and $LX\cap\ell'$. By the definition of the polar, $h$ maps each circle of positive radius of the pencil to itself. Since $h$ is a projective map, it sends chords of $\G$ tangent to $\g$ to chords of $\G$ tangent to $\g$, hence Poncelet polygons to Poncelet polygons, and $h(A_i)=A_i^*$. By \eqref{eq:LA} and the power of $L$,
\[
t_it_i^*=\frac{|LA_i|\cdot|LA_i^*|}{\mu^2}=\frac{(R^2-|OL|^2)\,d}{|OL|}.
\]
Since $|OL|$ is a root of $x^2-2x_0x+R^2=0$, we have $R^2-|OL|^2=2|OL|(x_0-|OL|)$, and the right-hand side equals $2d\sqrt{x_0^2-R^2}=\sqrt{(R^2+d^2-r^2)^2-4d^2R^2}=q$.

Finally let $n=2N$. On $\G$ the map $h$ is an involution without fixed points. It preserves both circles and the cyclic order of the points of $\G$, so it sends the tangent chord that goes from a point $A$ to $f(A)$ to the tangent chord that goes from $h(A)$ to $f(h(A))$; thus $hf=fh$. The maps $f$ and $h$ generate a finite abelian group of orientation preserving homeomorphisms of $\G$. Only the identity has fixed points in it, because an orientation preserving homeomorphism of a circle which has finite order and a fixed point is the identity. The cyclic order of an orbit shows that such a group is cyclic. A cyclic group has at most one element of order two, and $f^N$ is one. Hence $h=f^N$, that is, $A_i^*=A_{i+N}$, and the main diagonals pass through $L$.
\end{proof}

It follows that every invariant written in the $t_i$ stays an invariant, with the same value, when each $t_i$ is replaced by $q/t_i$.

\subsection{How the polygon moves}
\begin{lemma}\label{lem:kin}
Let a chord $XY$ of $\G$ move so that it stays tangent to a circle $\g'$ inside $\G$, at a point $C$ between $X$ and $Y$, with $\g'$ on the left of the chord directed from $X$ to $Y$. Let $\theta$ be the polar angle of $C$ seen from the center of $\g'$. Then, as $\theta$ increases, both endpoints move counterclockwise along $\G$, with speeds
\[
\Bigl|\frac{dX}{d\theta}\Bigr|=\frac{2R\,|XC|}{|XY|},\ \ \  \Bigl|\frac{dY}{d\theta}\Bigr|=\frac{2R\,|YC|}{|XY|}.
\]
\end{lemma}

\begin{proof}
Let $J$ and $\varrho$ be the center and the radius of $\g'$, and $u=(\cos\theta,\sin\theta)$, $v=(-\sin\theta,\cos\theta)$. Write $X=J+\varrho u+sv$ with $|s|=|XC|$. Then $dX/d\theta=(\varrho+s')v-su$, so the component of the velocity of $X$ normal to the chord has length $|s|$. The velocity is tangent to $\G$, and the angle between the chord and the tangent to $\G$ at $X$ has sine $|XY|/(2R)$. As for the direction, $s<0$ at $X$ and $s>0$ at $Y$, so the component $-s$ of the velocity along $u$ is positive at $X$ and negative at $Y$. The counterclockwise unit tangent to $\G$ has the same property, because the vector $v$ points into the disk at $X$ and out of it at $Y$.
\end{proof}

Applied to two consecutive sides of $\mathcal P$, which share the vertex $A_{i+1}$ with $|A_{i+1}B_i|=|A_{i+1}B_{i+1}|=t_{i+1}$ and move it in the same direction, the lemma gives
\begin{equation}\label{eq:poncelet-diff}
\frac{d\theta_i}{t_i+t_{i+1}}=\frac{d\theta_{i+1}}{t_{i+1}+t_{i+2}} .
\end{equation}
A relation of this kind is the starting point of Jacobi's proof of the porism \cite{BKOR,Jacobi}.

\subsection{What was known}
For the reader's convenience we list the results on Poncelet polygons between two circles that we are aware of and that are related to ours.
\begin{itemize}
\item Weill \cite{Weill1878}: the centroid of the contact points is fixed. Casey \cite{Casey} considers centroids of subsets of the contact points.
\item Weill \cite{Weill1880}: the ratio of the areas of the contact polygon and of $\mathcal P$ is constant.
\item Radi\'c \cite{Radic05,Radic14}: explicit relations between consecutive tangent lengths and many relations between $R$, $r$, $d$.
\item Radi\'c \cite{Radic04,Radic05}, Josefsson \cite{Josefsson,Josefsson15}, Cie\'slak and Mozgawa \cite{CM}: the polygons of largest and smallest area, or perimeter, in the family for $n=3,4,6,8$. Cie\'slak and Mozgawa conjectured the extremal polygons for every $n$ \cite[Section 7]{CM}.
\item Roitman, Garcia and Reznik \cite{RGR}: the sum of the cosines of the angles of $\mathcal P$ is constant; the pedal polygons of $\mathcal P$ with respect to $L$ and $L'$ have constant perimeter.
\item Keady and McAndrew \cite{KM}: for hexagons, $|A_1A_2||A_4A_5|+|A_2A_3||A_5A_6|+|A_3A_4||A_6A_1|$ is constant, together with related invariants built from sides and diagonals; the invariance of $e_2$, $e_4$, $e_6$, $e_3/e_1$, $e_5/e_1$ for the tangent lengths.
\item Hungerb\"uhler, Pohle and Zhang \cite{HPZ}: the excenters of $\mathcal P$ lie on a circle, and the ratio of the areas of $\mathcal P$ and of the polygon of excenters is constant.
\item Bellio, Garcia and Reznik \cite{BGR} and Garcia, Odehnal and Reznik \cite{GOR}: reciprocation with respect to $\G$, and families with several caustics from the pencil.
\end{itemize}
On the side of ellipses we shall meet the sum of cosines \cite{AST,BT}, the product of cosines of the tangential polygon \cite{BT}, the sum of the reciprocals of the focal distances \cite{RGH}, and the list of experimental invariants in \cite{RGK2}. The cyclic identities for Jacobi elliptic functions of Khare, Sukhatme and Lakshminarayan \cite{KS,KLS,KLS-local} are the analytic counterpart of several statements of this paper, since $t_i$ is, up to a factor, the function $\operatorname{dn}$ evaluated at equally spaced points.

\section{Proofs}\label{s:proofs}

\subsection{Weill's theorem for chains}
By a \emph{Poncelet chain} for a circle $\g'$ inside $\G$ we mean points $Q_1,\dots,Q_N$ of $\G$, with indices modulo $N$, such that every chord $Q_iQ_{i+1}$ touches $\g'$ at a point between its endpoints and $\g'$ lies on the left of every chord directed from $Q_i$ to $Q_{i+1}$. The chain may wind several times.

\begin{proposition}\label{prop:chain}
Let $\g'$ be a circle inside $\G$ whose center $J$ is different from $O$, and let a Poncelet chain $Q_1\dots Q_N$ for $\g'$ move continuously. Let $u_i$ be the tangent length from $Q_i$ to $\g'$. Then $\sum_iu_iu_{i+1}$ is constant.
\end{proposition}

\begin{proof}
Lemmas \ref{lem:tt} and \ref{lem:kin} did not use that the polygon is convex, so they apply to the chain, with $\delta=|OJ|$ and the radius $\varrho$ of $\g'$ in place of $d$ and $r$. Thus $u_iu_{i+1}=\text{const}+2\delta\varrho\cos\theta_i$, $u_i-u_{i+1}=2\delta\sin\theta_i$, and by \eqref{eq:poncelet-diff} there is a common differential $d\tau=d\theta_i/(u_i+u_{i+1})$. Therefore
\[
d\sum_i\cos\theta_i=-\sum_i\sin\theta_i\,(u_i+u_{i+1})\,d\tau=-\frac{d\tau}{2\delta}\sum_i\bigl(u_i^2-u_{i+1}^2\bigr)=0 .\qedhere
\]
\end{proof}

\subsection{Proof of Theorem \ref{thm:weill}}
Take $\g'=\g$ in Proposition \ref{prop:chain}: $W$ is constant. By \eqref{eq:tt} we have $\sum\sin\theta_i=0$ and $\sum\cos\theta_i=(W-nc)/(2dr)$, which gives (i), and $\sum|OB_i|^2=nR^2-W$, which is (ii). The signed distance from $X$ to the $i$-th side is $r-\overrightarrow{IX}\cdot(\cos\theta_i,\sin\theta_i)$; summing and using (i) gives (iii), and for $X=O$ the sum is $nr-d\sum\cos\theta_i$.

For (iv), multiply \eqref{eq:vieta} by $t_i$ and sum: $2W=2k\sum t_i^2/(t_i^2+r^2)=2k\,(n-r^2\sum|IA_i|^{-2})$. If $\widehat{A_i}$ is the angle of $\mathcal P$ at $A_i$, then $\sin(\widehat{A_i}/2)=r/|IA_i|$, so $\cos\widehat{A_i}=1-2r^2/|IA_i|^2$. For (v), the chord $B_{i-1}B_i$ of $\g$ subtends the angle $2\phi$ at $I$ with $\sin\phi=t_i/|IA_i|$, so $|B_{i-1}B_i|^2=4r^2t_i^2/(t_i^2+r^2)$.

For (vi), apply the invariance of $W$ to the conjugate polygon: $\sum1/(t_it_{i+1})=q^{-2}\sum t_i^*t_{i+1}^*=W/q^2$, because $\mathcal P^*$ belongs to the same family. Dividing \eqref{eq:p2} by $r^2p_i$ gives
\[
\frac{t_i}{t_{i+1}}+\frac{t_{i+1}}{t_i}=\frac{1}{r^2}\Bigl(2k-p_i-\frac{q^2}{p_i}\Bigr),
\]
and the sum over $i$ is $(2nk-2W)/r^2$. In each of (i) to (vi) the quantity on the left determines $W$ (in (iii) take $X=O$), so its constancy is equivalent to that of $W$. \qed

\subsection{Proof of Theorem \ref{thm:first}}
(a) The chord $A_iA_{i+2}$ subtends the inscribed angle $\widehat{A_{i+1}}$, so $|A_iA_{i+2}|=2R\sin\widehat{A_{i+1}}$. Since $\tan(\widehat{A_{i+1}}/2)=r/t_{i+1}$, we get $\sin\widehat{A_{i+1}}=2rt_{i+1}/(t_{i+1}^2+r^2)$, and by Lemma \ref{lem:vieta}
\[
|A_iA_{i+2}|=\frac{4Rr\,t_{i+1}}{t_{i+1}^2+r^2}=\frac{2Rr}{k}\,(t_i+t_{i+2}).
\]
Summing over $i$ and using \eqref{eq:S} gives \eqref{eq:diag2}. Dividing by the perimeter $2S/r$ gives $\Delta=2Rr/k$, which is \eqref{eq:euler-n}.

(b) The triangle $B_{i-1}IB_i$ has area $\tfrac12r^2\sin2\phi=r^3t_i/(t_i^2+r^2)$, which by \eqref{eq:vieta} equals $r^3(t_{i-1}+t_{i+1})/(2k)$. The sum over $i$ is $r^3\sum t_i/k=r^2S/k$.

(c) Sum \eqref{eq:diag-len} over $i$. For $n=2N$ and $m=N$ the sum over $i=1,\dots,N$ contains every $t_j$ once, which gives \eqref{eq:kvant}.

(d) Compute the area from $O$: the signed distance from $O$ to the $i$-th side is $r-d\cos\theta_i$, so $2S=\sum(t_i+t_{i+1})(r-d\cos\theta_i)$. Comparing with \eqref{eq:S} we get $\sum(t_i+t_{i+1})\cos\theta_i=0$, that is, $\sum(t_i+t_{i+1})(t_it_{i+1}-c)=0$. Now use $|OB_i|^2=R^2-t_it_{i+1}$ and $R^2-c=d^2+r^2$. \qed

\subsection{Proof of Theorem \ref{thm:second}}
Since $W_{n-m}=W_m$, we may assume that $m\le n/2$. Fix $m<n/2$ and let $u_i=\lambda_mt_i$ be the tangent lengths to $\g_m$. The points $A_i,A_{i+m},A_{i+2m},\dots$ form $\gcd(m,n)$ Poncelet chains for $\g_m$. Applying Proposition \ref{prop:chain} to each of them and adding, we see that $\sum_iu_iu_{i+m}=\lambda_m^2W_m$ is constant. If $n=2N$ then $t_it_{i+N}=t_it_i^*=q$ by Proposition \ref{prop:conj}. Since
\[
e_2=\sum_{1\le m<n/2}W_m+\begin{cases}\tfrac12W_{n/2},& n\text{ even},\\ 0,& n\text{ odd},\end{cases}
\]
$e_2$ is constant. Next, $S^2/r^2=(\sum t_i)^2=\sum t_i^2+2e_2$. Here $\sum t_i^2=\sum(|IA_i|^2-r^2)$; also $\sum|A_iA_{i+1}|^2=2\sum t_i^2+2W$; and $t_i^2=c+2\,\overrightarrow{IO}\cdot\overrightarrow{IA_i}$ because $|OA_i|=R$, so $\sum t_i^2=nc+2n\,\overrightarrow{IO}\cdot\overrightarrow{IG}$. This gives the three expressions. \qed

\subsection{Proof of Theorem \ref{thm:fourth}}
Summing \eqref{eq:p2} and using $\sum t_i^2=S^2/r^2-2e_2$,
\[
\sum p_i^2=-2r^2\sum t_i^2+2kW-nq^2=-2S^2+4r^2e_2+2kW-nq^2 .
\]
Since $|OB_i|^2=R^2-p_i$, we get $\sum|OB_i|^4=nR^4-2R^2W+\sum p_i^2$, which is \eqref{eq:four} because $k-R^2=-d^2$.

For the cubes multiply \eqref{eq:p2} by $p_i$:
\[
\sum p_i^3=-r^2\sum p_i(t_i^2+t_{i+1}^2)+2k\sum p_i^2-q^2W .
\]
Regrouping and using \eqref{eq:vieta},
\[
\sum p_i(t_i^2+t_{i+1}^2)=\sum t_i^3(t_{i-1}+t_{i+1})=2k\sum\frac{t_i^4}{t_i^2+r^2}=2k\Bigl(\sum t_i^2-nr^2+r^4\sum\frac1{|IA_i|^2}\Bigr).
\]
By Theorem \ref{thm:weill}(iv) the last sum is $(n-W/k)/r^2$. Collecting terms,
\[
\sum p_i^3=-6kr^2\sum t_i^2+(4k^2-q^2+2r^4)W-2knq^2,
\]
and $r^2\sum t_i^2=S^2-2r^2e_2$ gives the formula for $\sum p_i^3+6kS^2$. Finally
\[
\sum|OB_i|^6=nR^6-3R^4W+3R^2\sum p_i^2-\sum p_i^3=-6R^2S^2+6kS^2+\text{const}=-6d^2S^2+\text{const}. \qed
\]

\subsection{Proof of Theorem \ref{thm:product}}
Write $\dot t_i$ for the derivative with respect to $s$. By \eqref{eq:vieta},
\[
\dot W=\sum_i(t_{i-1}+t_{i+1})\,\dot t_i=k\sum_i\frac{2t_i\,\dot t_i}{t_i^2+r^2}=k\,\frac{\partial}{\partial s}\sum_i\log(t_i^2+r^2),
\]
and $t_i^2+r^2=|IA_i|^2$. Since $W$ is constant, so is $\prod|IA_i|$. The same computation for a Poncelet chain of $\g_m$, with Proposition \ref{prop:chain}, gives the statement about $I_m$. For even $n$ we have $\prod t_i=\prod_{i\le n/2}t_it_{i+n/2}=q^{n/2}$.

Let $n$ be odd and $m=(n-1)/2$. The points $A_i,A_{i+m},A_{i+2m},\dots$ form one chain for $\g_m$, with tangent lengths $u_i=\lambda_mt_i$, and Lemma \ref{lem:vieta} for this chain reads
\[
u_i+u_{i+2m}=\frac{2k_m\,u_{i+m}}{|I_mA_{i+m}|^2}.
\]
Since $2m\equiv-1$ modulo $n$, the left-hand side is $\lambda_m(t_i+t_{i-1})=\lambda_m|A_{i-1}A_i|$. Taking the product over $i$ and cancelling $\lambda_m^n$ gives \eqref{eq:odd-product}. \qed

\subsection{Proof of Theorem \ref{thm:limit}}
(a) The equalities are \eqref{eq:LA}. They give $|LA_i|:|LA_{i+1}|=t_i:t_{i+1}=|B_iA_i|:|B_iA_{i+1}|$, so $LB_i$ is the bisector of the angle at $L$ in the triangle $A_iLA_{i+1}$. The same argument works at $L'$.

(b) and (c). In the triangle $A_iLA_{i+m}$ the sides are $\mu t_i$, $\mu t_{i+m}$ and $\lambda_m(t_i+t_{i+m})$ by \eqref{eq:diag-len}. The law of cosines gives
\begin{equation}\label{eq:cosm}
\cos\angle A_iLA_{i+m}=\frac{\mu^2-\lambda_m^2}{2\mu^2}\Bigl(\frac{t_i}{t_{i+m}}+\frac{t_{i+m}}{t_i}\Bigr)-\frac{\lambda_m^2}{\mu^2}.
\end{equation}
For $m=1$ the sum over $i$ is given by Theorem \ref{thm:weill}(vi). For general $m$ we argue in the same way with $\g_m$ instead of $\g$: the tangent lengths $u_i=\lambda_mt_i$ of consecutive points of a chain satisfy a relation of the form \eqref{eq:p2}, the sum $\sum u_iu_{i+m}$ is constant by Proposition \ref{prop:chain}, and the sum $\sum1/(u_iu_{i+m})$ is constant because the homology $h$ of Proposition \ref{prop:conj} preserves $\g_m$ and maps chains to chains. At $L'$ replace $\mu$ by $\mu'$.

(d) By \eqref{eq:cosm} with $m=1$ and the formula in the proof of Theorem \ref{thm:weill}(vi), $\cos\varphi_i=a-b\,(p_i+q^2/p_i)$ with $b=(\mu^2-1)/(2\mu^2r^2)$ and a constant $a$. Hence
\[
\sum\cos^2\varphi_i=\text{const}+b^2\Bigl(\sum p_i^2+q^4\sum\frac1{p_i^2}\Bigr).
\]
By Theorem \ref{thm:fourth}, $\sum p_i^2=\text{const}-2S^2$. The numbers $q/t_i$ are the tangent lengths of $\mathcal P^*$, so $q^4\sum p_i^{-2}=\sum(p_i^*)^2=\text{const}-2S^{*2}$ with the same constant. For even $n$ the polygon $\mathcal P^*$ is $\mathcal P$ with shifted labels.

(e) By (a), at $L$ the ray $LB_i$ bisects the angle $\varphi_i$, hence $\psi_i=(\varphi_i+\varphi_{i+1})/2$ and $\cos2\psi_i=\cos\angle A_iLA_{i+2}$. At $L'$ all rays $L'A_i$ lie in a half-plane, the direction of $L'B_i$ is the mean of the directions of $L'A_i$ and $L'A_{i+1}$, and again $\cos2\psi_i'=\cos\angle A_iL'A_{i+2}$. So we have $\sum\cos^2\psi_i=n/2+\tfrac12\sum\cos\angle A_iLA_{i+2}$, which is constant by (c) if $n\ge5$, by (b) if $n=3$, and trivially at $L$ if $n=4$ since then $\psi_i=\pi/2$.

For the product use \eqref{eq:cosm} with $m=2$, which is valid for every $n$ with $\lambda_2=2Rr/k$ by the proof of Theorem \ref{thm:first}(a):
\[
\cos^2\psi_i=\frac{1+\cos\angle A_iLA_{i+2}}{2}=\frac{(\mu^2-\lambda_2^2)(t_i+t_{i+2})^2}{4\mu^2\,t_it_{i+2}} .
\]
By \eqref{eq:vieta}, $\prod(t_i+t_{i+2})=(2k)^n\prod t_i/\prod|IA_i|^2$. Taking the product over $i$ gives the formula in (e). \qed

\section{Doing without calculus}\label{s:nocalc}

Derivatives were used in two places: in Lemma \ref{lem:kin}, for the proof of Proposition \ref{prop:chain}, and in the identity \eqref{eq:dW} of Theorem \ref{thm:product}. Here is another proof of Proposition \ref{prop:chain} in which two positions of the polygon are compared directly. It rests on a lemma which is essentially Poncelet's own main lemma \cite{BKOR,GOR}, and which we prove in order to keep track of the constants.

\begin{lemma}[Four points of contact]\label{lem:four}
Let $A$, $A'$, $A_1$, $A_1'$ lie on $\G$ in this counterclockwise order, and suppose the chords $AA_1$ and $A'A_1'$ touch $\g$ at $B$ and $B'$, with $\g$ on the left of both chords when they are directed from $A$ to $A_1$ and from $A'$ to $A_1'$. Let the line $BB'$ meet $AA'$ at $C$ and $A_1A_1'$ at $C'$. Then:
\begin{enumerate}[label=\textup{(\arabic*)}]
\item the points lie on the line in the order $C$, $B$, $B'$, $C'$;
\item there is a number $\lambda$ such that $|AC|=\lambda|AB|$, $|A'C|=\lambda|A'B'|$, $|A_1C'|=\lambda|A_1B|$ and $|A_1'C'|=\lambda|A_1'B'|$;
\item the circle $\g'$ of the pencil of $\G$ and $\g$ whose center $J$ satisfies $\overline{OJ}=\lambda^2\,\overline{OI}$ touches $AA'$ at $C$ and $A_1A_1'$ at $C'$;
\item if $\varrho$ is the radius of $\g'$, then
\begin{equation}\label{eq:four-powers}
\pw_\G(B)-\pw_\G(B')=-\frac{r}{\lambda\varrho}\,\bigl(\pw_\G(C')-\pw_\G(C)\bigr).
\end{equation}
\end{enumerate}
\end{lemma}

\begin{figure}[ht]
\centering
\begin{tikzpicture}[scale=0.85]
\draw (0.9,0) circle (3.0);
\draw (0,0) circle (1.5);
\draw[dashed] (0.625,0) circle (2.601);
\draw[thick] (3.533,1.438)--(-1.681,1.530);
\draw[thick] (0.513,2.975)--(-2.095,-0.179);
\draw (3.533,1.438)--(0.513,2.975);
\draw (-1.681,1.530)--(-2.095,-0.179);
\draw[thin] (2.5,2.638)--(-2.6,0.291);
\foreach \p in {(0.026,1.500),(-1.156,0.956),(1.805,2.318),(-1.903,0.612)} \fill \p circle (1.6pt);
\fill (0.9,0) circle (1.4pt) node[below] {$O$};
\fill (0,0) circle (1.4pt) node[below] {$I$};
\fill (0.625,0) circle (1.4pt) node[above] {$J$};
\node[right] at (3.533,1.438) {$A$};
\node[above] at (0.513,2.975) {$A'$};
\node[above left] at (-1.681,1.530) {$A_1$};
\node[left] at (-2.095,-0.179) {$A_1'$};
\node[below] at (0.15,1.48) {$B$};
\node[below right] at (-1.156,0.956) {$B'$};
\node[above right] at (1.805,2.318) {$C$};
\node[left] at (-1.95,0.72) {$C'$};
\node at (3.4,-2.1) {$\G$};
\node at (0.9,-0.85) {$\g$};
\node at (2.55,-1.3) {$\g'$};
\end{tikzpicture}
\caption{Lemma \ref{lem:four}: the four points of contact are collinear.}\label{fig:two}
\end{figure}
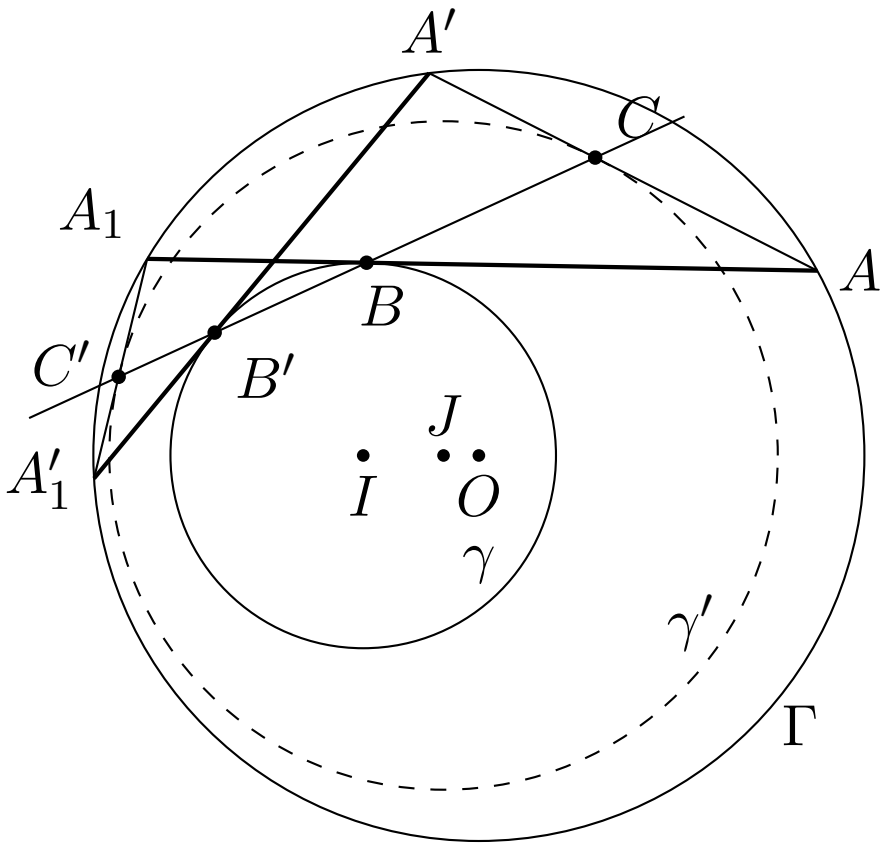

\begin{proof}
The endpoints of the two chords alternate on $\G$, so the chords meet at a point $Q$. The circle $\g$ lies to the left of both directed chords, hence in the angle $AQA_1'$, so $B$ lies on the segment $QA$ and $B'$ on the segment $QA_1'$, and $|QB|=|QB'|$. Let $\eta$ be the base angle of the isosceles triangle $QBB'$. The line $BB'$ is parallel to the bisector of the angle $AQA'$; beyond $B$ it therefore leaves the triangle $QAA'$ through the side $AA'$, and beyond $B'$ it leaves the triangle $QA_1'A_1$ through $A_1'A_1$. This proves (1).

In each of the four triangles $ABC$, $A'B'C$, $A_1BC'$, $A_1'B'C'$ the angle at $B$ or $B'$ equals $\eta$ or is vertical to an angle $\eta$. The angles of $ABC$ and $A'B'C$ at $C$ are supplementary. By the law of sines
\[
\frac{|AC|}{|AB|}=\frac{\sin\eta}{\sin\angle ACB}=\frac{|A'C|}{|A'B'|}.
\]
The angles $BAC$ and $B'A_1'C'$ are inscribed in $\G$ and subtend the same arc $A'A_1$, so the triangles $ABC$ and $A_1'B'C'$ have two equal angles and $\angle ACB=\angle A_1'C'B'$. The angles of $A_1BC'$ and $A_1'B'C'$ at $C'$ are supplementary. Hence the same ratio $\lambda=\sin\eta/\sin\zeta$, where $\zeta=\angle ACB$, serves for all four triangles. This is (2).

Let $\g'$ be as in (3); a priori it might have no real points, but its power function is still defined by \eqref{eq:pencil}. By Lemma \ref{lem:casey} and (2), $\pw_{\g'}(A)=\lambda^2|AB|^2=|AC|^2$ and $\pw_{\g'}(A')=|A'C|^2$. On the line $AA'$, with a coordinate $x$, the power with respect to $\g'$ is a quadratic polynomial with leading coefficient $1$. It agrees with $(x-x_C)^2$ at the two points $A$ and $A'$, and so it is equal to $(x-x_C)^2$. Thus $\g'$ is a real circle touching $AA'$ at $C$. The same argument on the line $A_1A_1'$ shows that $\g'$ touches it at $C'$.

(4) Introduce a coordinate on the line $BB'$ and let $b<b'$ and $c<c'$ be the coordinates of the four points; then $c<b<b'<c'$. For a point $X$ of the line with coordinate $x$ we have $\pw_\g(X)=(x-b)(x-b')$ and $\pw_{\g'}(X)=(x-c)(x-c')$. By \eqref{eq:three-powers},
\[
(\lambda^2-1)\,\pw_\G(X)=\lambda^2(x-b)(x-b')-(x-c)(x-c').
\]
Put $D_B=\pw_\G(B)-\pw_\G(B')$ and $D_C=\pw_\G(C')-\pw_\G(C)$. Evaluating at the four points we find
\[
(\lambda^2-1)\,D_B=(b'-b)(b+b'-c-c'),\ \ \ (\lambda^2-1)\,D_C=-\lambda^2(c'-c)(b+b'-c-c').
\]
Multiply the first equation by $\lambda^2(c'-c)$, the second by $b'-b$, and add:
\[
(\lambda^2-1)\bigl(\lambda^2(c'-c)\,D_B+(b'-b)\,D_C\bigr)=0 .
\]
Here $\lambda^2\ne1$, since otherwise $\g'=\g$ and the line $BB'$ would meet $\g$ at three points $C$, $B$, $B'$. Also $\varrho>0$, because $C\ne C'$. The chord $BB'$ of $\g$ makes the angle $\eta$ with the tangent $AA_1$, and the chord $CC'$ of $\g'$ makes the angle $\zeta$ with the tangent $AA'$. So $b'-b=2r\sin\eta$, $c'-c=2\varrho\sin\zeta$, and $(b'-b)/(c'-c)=\lambda r/\varrho$. This gives \eqref{eq:four-powers}.
\end{proof}

\begin{proposition}\label{prop:finite}
Let $\mathcal P=A_1\dots A_n$ and $\mathcal P'=A_1'\dots A_n'$ be two Poncelet polygons for $\G$ and $\g$, with points of contact $B_i$ and $B_i'$. Then $\sum|OB_i|^2=\sum|OB_i'|^2$.
\end{proposition}

\begin{proof}
We may assume that the polygons are different. Label the second one so that $A_1'$ lies on the open arc $A_1A_2$; since the Poncelet map preserves the cyclic order of the points of $\G$, each $A_i'$ then lies on the open arc $A_iA_{i+1}$. For each $i$ apply Lemma \ref{lem:four} to the chords $A_iA_{i+1}$ and $A_i'A_{i+1}'$. We get a number $\lambda_{(i)}$, a circle $\g'_{(i)}$ of the pencil, and points $C_i$ on $A_iA_i'$ and $C_i'$ on $A_{i+1}A_{i+1}'$.

By part (2), the point $C_i'$ divides the segment $A_{i+1}A_{i+1}'$ in the ratio $t_{i+1}:t_{i+1}'$, and so does $C_{i+1}$. Hence $C_i'=C_{i+1}$, and then $\lambda_{(i)}=|A_{i+1}C_{i+1}|/t_{i+1}=\lambda_{(i+1)}$. So all the $\lambda_{(i)}$ are equal to one number $\lambda$, all the circles $\g'_{(i)}$ coincide with one circle $\g'$, which touches every line $A_iA_i'$ at $C_i$, and \eqref{eq:four-powers} reads
\[
\pw_\G(B_i)-\pw_\G(B_i')=-\frac{r}{\lambda\varrho}\,\bigl(\pw_\G(C_{i+1})-\pw_\G(C_i)\bigr).
\]
The sum of the right-hand sides over $i$ is zero. Since $\pw_\G(X)=|OX|^2-R^2$, the proposition follows.
\end{proof}

Proposition \ref{prop:finite} is Weill's theorem in the form (ii) of Theorem \ref{thm:weill}. The same argument proves Proposition \ref{prop:chain} for an arbitrary Poncelet chain and its circle. If the chain winds several times, one compares two positions that are close to each other, so that the endpoints of corresponding chords alternate as required in Lemma \ref{lem:four}. This shows that $\sum u_iu_{i+1}$ is locally constant, and it is constant because the family is connected. Along the way it shows that the lines $A_iA_i'$ touch one circle of the pencil, which is a case of Poncelet's general theorem.

With this, Theorems \ref{thm:first}, \ref{thm:weill}, \ref{thm:second}, \ref{thm:fourth} and \ref{thm:limit}, except for the product in its part (e), are proved without derivatives or integrals, given Theorem \ref{thm:PGT}. The exception is Theorem \ref{thm:product}, whose proof is the identity \eqref{eq:dW}; Corollary \ref{cor:more}(a) gives another proof that $\prod|IA_i|$ is constant, by complex analysis.

\begin{question}
Is there a proof that $\prod|IA_i|$ is constant in the spirit of Proposition \ref{prop:finite}?
\end{question}

\section{Confocal conics}\label{s:conics}

\subsection{The dictionary}
Let $\omega$ be a circle with center $L$ and radius $\rho$, and let us reciprocate in $\omega$ as in Section \ref{ss:intro-polar}. The focus $F_1=L$ is common to the ellipses $E=\g^\circ$ and $E'=\G^\circ$. Let $F_2$ be their second focus, $a>b$ the semi-axes of $E$ and $a'>b'$ those of $E'$. One checks that
\[
\frac{b^2}{a}=\frac{\rho^2}{r},\ \ \ \frac{b'^2}{a'}=\frac{\rho^2}{R},\ \ \  \mu^2=\frac{|OL|}{|OI|}=\frac{b^2}{b^2-b'^2}.
\]
The polar of the vertex $A_i$ is a line $s_i$ tangent to $E'$, the pole of the side $A_iA_{i+1}$ is a point $P_i$ of $E$, and $P_i$ lies on $s_i$ and $s_{i+1}$. So $P_1P_2\dots P_n$ is a polygon inscribed in $E$ and circumscribed about $E'$. The tangent to $E$ at $P_i$ is the polar of $B_i$. Since $s_i\perp LA_i$, $s_{i+1}\perp LA_{i+1}$ and the tangent is perpendicular to $LB_i$, Theorem \ref{thm:limit}(a) says that the two sides at $P_i$ make equal angles with the tangent to $E$. So $P_1\dots P_n$ is a billiard trajectory in $E$. It is a convex polygon, that is, a trajectory of winding number one, and every periodic billiard trajectory of winding number one arises in this way \cite{RGR}.

Let $h_i$ and $g_i$ be the distances from $F_1$ and $F_2$ to the line $s_i$. Then
\begin{equation}\label{eq:dict}
h_i=\frac{\rho^2}{|LA_i|}=\frac{\rho^2}{\mu\,t_i},\ \ \  g_i=\frac{b'^2}{h_i}=\frac{b'^2\mu}{\rho^2}\,t_i ,
\end{equation}
because the product of the distances from the foci to a tangent of $E'$ is $b'^2$. Thus \emph{the tangent lengths of the bicentric polygon are proportional to the distances from the second focus to the lines of the sides of the billiard trajectory}. Here and below the distances are taken to the lines that contain the sides; the foot of a perpendicular may lie outside the side. The conjugation $t_i\mapsto q/t_i$ of Proposition \ref{prop:conj} becomes the central symmetry of the ellipse, which interchanges the foci.

\begin{lemma}\label{lem:focal}
Let $P$ be a point outside an ellipse with foci $F_1$, $F_2$ and minor semi-axis $b'$, and let $g$, $g^+$ be the distances from $F_2$ to the two tangents from $P$. Then $g\,g^+=b'^2\,|PF_2|/|PF_1|$.
\end{lemma}

\begin{proof}
Let $h$, $h^+$ be the distances from $F_1$ to the same tangents. The two tangents from $P$ are isogonal with respect to the angle $F_1PF_2$, so $g/h^+=|PF_2|/|PF_1|=g^+/h$. Multiply these equalities and use $gh=g^+h^+=b'^2$.
\end{proof}

By the lemma and \eqref{eq:dict}, the product $t_it_{i+1}=R^2-|OB_i|^2$ corresponds to the ratio of the focal distances of $P_i$:
\begin{equation}\label{eq:dict2}
g_ig_{i+1}=b'^2\,\frac{|F_2P_i|}{|F_1P_i|}=b'^2\Bigl(\frac{2a}{|F_1P_i|}-1\Bigr),\ \ \  \sum_i g_i=\frac{b'^2\mu}{\rho^2r}\,S .
\end{equation}

\subsection{Billiards in ellipses}
Let $P_1\dots P_n$ be a periodic billiard trajectory of winding number one in the ellipse $E$, with the confocal ellipse $E'$ as caustic, let $g_i$ be the distance from the focus $F_2$ to the line $P_{i-1}P_i$, let $\beta_i$ be the interior angle of the trajectory at $P_i$, and let $\beta_i^{\,\rm t}$ be the angles of the polygon formed by the tangents to $E$ at the points $P_i$.

\begin{theorem}\label{thm:ellipse}
Over the family of periodic billiard trajectories of winding number one with a given caustic $E'$, which is an ellipse, the following are constant:
\begin{enumerate}[label=\textup{(\alph*)}]
\item $\sum1/|F_1P_i|$ \cite{RGH,RGR};
\item $\sum\cos\beta_i$ \cite{AST,BT};
\item $\prod\cos\beta_i^{\,\rm t}$ \cite{BT};
\item $\sum_ig_ig_{i+m}$ for $1\le m\le n-1$, and hence $\bigl(\sum g_i\bigr)^2-\sum g_i^2$;
\item $\displaystyle\sum_{i=1}^n\frac1{|F_1P_i|^2}+\frac{1}{2b^2(b^2-b'^2)}\Bigl(\sum_{i=1}^ng_i\Bigr)^2$.
\end{enumerate}
\end{theorem}

\begin{proof}
(a) By \eqref{eq:dict2}, $\sum1/|F_1P_i|$ is an affine function of $\sum g_ig_{i+1}$, which is a multiple of $W$. (b) The outer normals of the sides $s_i$ and $s_{i+1}$ of the convex polygon $P_1\dots P_n$ have the directions of $\overrightarrow{LA_i}$ and $\overrightarrow{LA_{i+1}}$. The angle between them is $\varphi_i$, so $\beta_i=\pi-\varphi_i$, and we apply Theorem \ref{thm:limit}(b). (c) The tangents to $E$ at $P_i$ and $P_{i+1}$ are perpendicular to $LB_i$ and $LB_{i+1}$, so $|\cos\beta_i^{\,\rm t}|=|\cos\psi_i|$, and we apply Theorem \ref{thm:limit}(e), with the remark on the signs made there; if the product vanishes for one trajectory, it vanishes for all. (d) is Theorem \ref{thm:second}. (e) By Theorem \ref{thm:fourth} and \eqref{eq:S}, $\sum p_i^2+2r^2(\sum t_i)^2$ is constant. Substituting \eqref{eq:dict} and \eqref{eq:dict2}, dividing by $4a^2b'^4$ and discarding the term with $\sum1/|F_1P_i|$, which is constant by (a), we get (e) after a short computation with $\mu^2=b^2/(b^2-b'^2)$ and $\rho^2=rb^2/a$.
\end{proof}

Parts (a), (b), (c) are known theorems. What is perhaps new is that all three come from the theorem of Weill: (a) and Weill's theorem are literally the same statement, (b) needs in addition the central symmetry of the ellipse, and (c) needs the identity \eqref{eq:dW}. Together with Section \ref{s:nocalc} this gives a proof of (a) and (b) that uses no calculus and no elliptic functions.

Part (d) has a pleasant form in terms of pedal polygons. The feet $Q_i$ of the perpendiculars from $F_2$ to the lines of the sides lie on the auxiliary circle of the caustic, with center at the center of the ellipse and radius $a'$, and $g_i=|F_2Q_i|$. So for the polygon $Q_1\dots Q_n$ inscribed in this circle the sums $\sum_i|F_2Q_i|\cdot|F_2Q_{i+m}|$, $1\le m\le n-1$, are constant. The inversion in $\omega$ sends $\G$ to this circle and $A_i$ to the foot of the perpendicular from $F_1$ to $s_i$; the feet of the perpendiculars from $F_2$ are obtained from these by the central symmetry of the ellipse, which preserves the family of trajectories.

\subsection{Hyperbolas}
If we reciprocate in a circle of radius $\rho$ centered at the outer limiting point $L'$, the circles $\G$ and $\g$ become confocal hyperbolas $H'$ and $H$ with a focus at $L'$. The polar of $A_i$ is a line $s_i$ tangent to $H'$, and the pole of the side $A_iA_{i+1}$ is the common point $P_i$ of $s_i$ and $s_{i+1}$, a point of $H$. Some care is needed here. The two tangents from $L'$ to $\g$ meet $\G$, so they are sides of polygons of the family, and when the side $A_iA_{i+1}$ passes through $L'$, the point $P_i$ is at infinity and the lines $s_i$, $s_{i+1}$ are parallel. So the reciprocal polygon lives in the projective plane. Its vertices may lie on both branches of $H$, and it is a billiard trajectory only in a generalized sense, see \cite{RGR}.

The following quantities, however, are defined for every position of the bicentric polygon: the distance $h_i=\rho^2/|L'A_i|=\rho^2/(\mu't_i)$ from the focus $L'$ to the line $s_i$; the unit vector $u_i=\overrightarrow{L'A_i}/|L'A_i|$, which is the unit normal of $s_i$ pointing from the focus towards this line; and the unit vector $w_i=\overrightarrow{L'B_i}/|L'B_i|$, which for a finite vertex $P_i$ is the unit normal of the tangent to $H$ at $P_i$, pointing from the focus towards this tangent. In these terms we get, for instance:
\begin{itemize}
\item $\sum_iu_i\cdot u_{i+1}$ is constant, by Theorem \ref{thm:limit}(b). The angle between $u_i$ and $u_{i+1}$ may be obtuse, and the statement fails if one takes instead the acute angles between the lines $s_i$ and $s_{i+1}$;
\item $\sum1/(h_ih_{i+m})$ and $\sum h_ih_{i+m}$ are constant for $1\le m\le n-1$, by Theorem \ref{thm:second} and Proposition \ref{prop:conj};
\item $\prod_iw_i\cdot w_{i+1}$ is constant, by Theorem \ref{thm:limit}(e).
\end{itemize}
Statements that involve sums of angles, or distances with a sign, need more care at $L'$ than at $L$, because $L'$ lies outside the polygon.

\subsection{Which side is nicer?}
A few comments on taste. Weill's theorem looks best between two circles: a centroid is fixed. Its billiard form (a) is also simple. The problem with the fourth powers, on the other hand, clearly prefers circles. Between circles it needs only the points $B_i$, the center $O$ and the area, and the coefficient is the number $2$. On the ellipse it becomes (e), which mixes the two foci and has a coefficient depending on both conics. The second quadratic invariant is acceptable on both sides, as $S^2-r^2\sum|IA_i|^2$ and as (d). The sum of cosines is the opposite case: it is natural for billiards, and between circles it turns into the less obvious Theorem \ref{thm:limit}(b) about the angles at a limiting point. We have not found this last statement in the literature in this form, though it is implicit in \cite{AST,RGR}.

\section{How the statements were found}\label{s:found}

\subsection{Experiment}
The identities of this paper were first found numerically. One fixes $R$ and $d$, finds $r$ for which the polygon closes after $n$ steps, moves the first vertex, and watches which combinations of lengths, areas and angles stay put to twelve digits. The two circles are convenient for this, since a polygon is constructed by intersecting lines and circles. A warning is in order. Some quantities that are \emph{not} constant, like $\sum\cos^2\angle A_iLA_{i+1}$, varied only in the fifth or sixth digit in several of the families with odd $n$ that we tried, and a picture on the screen will not show this.

Once a candidate is found, it helps to translate it. We went back and forth between four languages: lengths and areas between the two circles, the tangent lengths $t_i$, angles at the limiting points, and the confocal ellipses. A statement that is opaque in one language can be a known theorem in another, as with Theorem \ref{thm:limit}(b) and the sum of cosines for billiards, or it can become easy, as with the Kvant problem \eqref{eq:kvant}, which in the variables $t_i$ is the remark that $|A_iA_{i+N}|=\mu(t_i+t_{i+N})$.

\subsection{The elliptic curve}\label{ss:curve}
There is a way to predict what will be constant. The pairs $(t_i,t_{i+1})$ satisfy \eqref{eq:Phi}, a symmetric equation $\Phi(x,y)=0$ of degree two in each variable. Over the complex numbers it defines a curve $\mathcal E$ of genus one in the product of two projective lines, the map $(t_{i-1},t_i)\mapsto(t_i,t_{i+1})$ is a translation $T$ of $\mathcal E$ with $T^n=\mathrm{id}$, and any rational expression in the $t_i$ that is invariant under cyclic shifts is a meromorphic function on $\mathcal E$. Such a function is constant as soon as it has no poles. This is the method of \cite{AST,RGR}, and it goes back to Jacobi.

In more detail, write \eqref{eq:Phi} as
\begin{equation}\label{eq:curve}
(x^2+r^2)\,y^2-2kxy+r^2x^2+q^2=0 .
\end{equation}
The discriminant with respect to $y$ is $4r^2(b_0^2-x^2)(x^2-a_0^2)$, where $a_0^2=(R-d)^2-r^2$ and $b_0^2=(R+d)^2-r^2$. So there are four distinct branch points, and $\mathcal E$ is a smooth curve of genus one. By Lemma \ref{lem:vieta}, $t_{i-1}$ and $t_{i+1}$ are the two roots of \eqref{eq:curve} with $x=t_i$. Thus $T$ is the composition of the involution that exchanges the two roots $y$ and the involution that exchanges $x$ and $y$. Both have fixed points, so $T$ is a translation. By Lemma \ref{lem:tt} a pair $(t_i,t_{i+1})$ determines $\theta_i$, hence the side and the vertex $A_i$, so $T$ has order exactly $n$. The function $t_0=x$ has exactly two poles $Q_+$ and $Q_-$ on $\mathcal E$, both simple, at which $t_1=ir$ and $t_1=-ir$.

\begin{proposition}\label{prop:liouville}
Let $F(t_1,\dots,t_n)$ be a polynomial which is invariant under cyclic shifts, under reversal of the indices and under the change of sign of all the variables, and in which no variable occurs to a power higher than one. Then $F$ is constant over the Poncelet family.
\end{proposition}

The main examples are homogeneous polynomials of even degree with the first two symmetries. Without the third symmetry the statement is false: $e_2+e_1$ is not constant.

\begin{proof}
Let $\iota$ be the map $(x,y)\mapsto(-x,-y)$ of $\mathcal E$. It commutes with $T$ and exchanges $Q_+$ and $Q_-$. We claim that $t_j$ is finite at $Q_+$ if $j\not\equiv0$ modulo $n$. Indeed, $t_j(Q_+)=\infty$ means that $T^jQ_+$ is $Q_+$ or $Q_-$. The first is impossible because $T$ is a translation of order $n$. If $T^jQ_+=Q_-=\iota Q_+$, then $T^{2j}Q_+=Q_+$, so $n$ is even and $j\equiv n/2$. But for even $n$ the relation $t_{n/2}=q/t_0$ of Proposition \ref{prop:conj} holds for the real polygons, hence on $\mathcal E$, and so $t_{n/2}$ vanishes at $Q_+$.

By Lemma \ref{lem:vieta}, $t_{-1}=-t_1=-ir$ at $Q_+$. The sequence $t_j'=-t_{-j}$ is again an orbit of $T$, and at $Q_+$ it has the same initial pair $(t_0',t_1')=(\infty,ir)$ as the sequence $(t_j)$. Hence $t_{-j}=-t_j$ at $Q_+$ for all $j\not\equiv0$.

Now write $F=t_0\cdot C+D$, where $C$ and $D$ are polynomials in the $t_j$ with $j\not\equiv0$; they are holomorphic at $Q_+$. The symmetries of $F$ show that $C$ changes sign when all the variables change sign, and that $C$ is invariant under $j\mapsto-j$. So the relation $t_{-j}=-t_j$ gives $C(Q_+)=-C(Q_+)=0$, and the simple pole of $t_0$ is cancelled. Thus $F$ is holomorphic at $Q_+$, and also at $Q_-$ because $F\circ\iota=F$. All the other poles of the functions $t_j$ are at the points $T^{-j}Q_\pm$, and $F\circ T=F$. So $F$ has no poles.
\end{proof}

Proposition \ref{prop:liouville} contains Theorems \ref{thm:weill} and \ref{thm:second}, and it was our guide to them. It also gives invariants for which we have no elementary proof.

\begin{corollary}\label{cor:more}
The following are constant over the Poncelet family:
\begin{enumerate}[label=\textup{(\alph*)}]
\item every elementary symmetric function $e_{2j}(t_1,\dots,t_n)$ of even degree;
\item $\sum_i|OB_i|^2\,|OB_{i+j}|^2$ for $2\le j\le n-2$, and more generally the sum, over all rotations and reflections of a set $J$ of pairwise non-adjacent sides, of $\prod_{i\in J}|OB_i|^2$;
\item $\sum_i|OB_i|^2\,|OB_{i+1}|^2-S^2$;
\item for odd $n$, the ratio $S/(t_1t_2\cdots t_n)$.
\end{enumerate}
\end{corollary}

Indeed, (a) is immediate; (b) follows since $R^2-|OB_i|^2=t_it_{i+1}$ and the sides are non-adjacent; (c) follows from (b), Theorem \ref{thm:weill}(ii) and \eqref{eq:four} by expanding $(\sum|OB_i|^2)^2$; and (d) is (a) for the conjugate polygon, because $e_1(t)/e_n(t)=e_{n-1}(1/t)$. Note that the constancy of $\prod|IA_i|$ also follows from (a), since
\[
\prod_{j=1}^n\Bigl(1+\frac{i\,t_j}{r}\Bigr)=-\frac{1}{r^n}\prod_{j=1}^n|IA_j|
\]
is real and its real part is $\sum_j(-1)^je_{2j}/r^{2j}$. Here we used that the argument of $1+i\,t_j/r$ is the angle $B_{j-1}IA_j$, and that these angles add up to $\pi$.

The same reasoning explains the shape of Theorem \ref{thm:fourth}. Let $\mathcal E_1$ be the quotient of $\mathcal E$ by the shifts. On $\mathcal E_1$ the function $S=r\sum t_j$ has two simple poles, the images of $Q_+$ and $Q_-$, so it has degree two. The reversal of the indices preserves $S$ and exchanges these two poles. Hence the quotient of $\mathcal E_1$ by the reversal is a projective line with coordinate $S$. All the poles of a polynomial in the $t_i$ lie on the orbits of $Q_+$ and $Q_-$, which go to the point $S=\infty$ of this line. It follows that a polynomial in the $t_i$ which is invariant under shifts and reversal is a polynomial in $S$, and if it is also invariant under the change of sign of all the $t_i$, it is a polynomial in $S^2$, of degree equal to half the order of its pole at $Q_+$. For $\sum|OB_i|^4$ and $\sum|OB_i|^6$ the pole has order two, so they are affine functions of $S^2$. For $\sum|OB_i|^8$ the pole has order four. Its leading term at $Q_+$ comes from $p_{-1}^4+p_0^4=t_0^4(t_{-1}^4+t_1^4)$, which behaves like $2r^4t_0^4$, while $S^4$ behaves like $r^4t_0^4$. So
\[
\sum_{i=1}^n|OB_i|^8=2S^4+\beta S^2+\text{const},
\]
where $\beta$ depends on $n$ and on the circles; it is computed in Proposition \ref{prop:ten} below.

The same quotient gives all the elementary symmetric functions of the tangent lengths at once. Let $r_m$ be the radius of $\g_m$, and let $\varepsilon=0$ for odd $n$ and $\varepsilon=1$ for even $n$.

\begin{theorem}\label{thm:poly}
There are polynomials $\mathcal A(z)$ and $\mathcal B(z)$, depending only on $n$ and on the two circles, such that for every polygon of the family
\begin{equation}\label{eq:poly}
\prod_{j=1}^n(z+t_j)=\mathcal A(z)+\frac{S}{r}\,\mathcal B(z).
\end{equation}
In other words, every $e_{2j}(t_1,\dots,t_n)$ and every ratio $e_{2j+1}(t_1,\dots,t_n)/e_1(t_1,\dots,t_n)$ is constant. Moreover,
\begin{equation}\label{eq:B}
\mathcal B(z)=z^{\varepsilon}\prod_{1\le m<n/2}\Bigl(z^2+\frac{r_m^2}{\lambda_m^2}\Bigr),\ \ \ \text{in particular}\ \ \ \frac{e_3}{e_1}=\sum_{1\le m<n/2}\frac{r_m^2}{\lambda_m^2}.
\end{equation}
\end{theorem}

\begin{proof}
In the polynomial $e_j$ every variable occurs to the first power, so its poles at $Q_+$ and $Q_-$ are at most simple, like those of $S$. It is invariant under shifts and reversal. Hence, by what was said above, $e_j$ is a polynomial of degree at most one in $S$. The change of sign of all the $t_i$ multiplies $e_j$ by $(-1)^j$ and $S$ by $-1$. So $e_j$ is constant for even $j$ and a constant multiple of $S$ for odd $j$, and this is \eqref{eq:poly}.

Let $\tau_j$ be the value of $t_j$ at $Q_+$ for $j\not\equiv0$. Since $S/(r\,t_0)$ tends to $1$ at $Q_+$, dividing \eqref{eq:poly} by $t_0$ and passing to $Q_+$ we get $\mathcal B(z)=\prod_{j\not\equiv0}(z+\tau_j)$. We saw in the proof of Proposition \ref{prop:liouville} that $\tau_{n-j}=-\tau_j$ and that $\tau_{n/2}=0$ for even $n$. Thus $\mathcal B(z)=z^{\varepsilon}\prod_{1\le m<n/2}(z^2-\tau_m^2)$. To find $\tau_m$, write the relation \eqref{eq:curve} for a chain of $\g_m$, whose tangent lengths are $\lambda_mt_i$, and divide it by $\lambda_m^4$:
\[
t_i^2t_{i+m}^2+\frac{r_m^2}{\lambda_m^2}\,(t_i^2+t_{i+m}^2)-\frac{2k_m}{\lambda_m^2}\,t_it_{i+m}+q^2=0,\ \ \ k_m=R^2-|OI_m|^2 .
\]
The last term is $q^2$, because the product of the smallest and the largest tangent length from the points of $\G$ to $\g_m$ is $\lambda_m^2q$. The relation holds on $\mathcal E$. Dividing by $t_0^2$ and passing to $Q_+$ we get
\begin{equation}\label{eq:taum}
\tau_m^2=-\frac{r_m^2}{\lambda_m^2},
\end{equation}
which gives the formula for $\mathcal B$. The formula for $e_3/e_1$ is the coefficient of $z^{n-3}$ in \eqref{eq:poly}.
\end{proof}

The constancy of the $e_{2j}$ is Corollary \ref{cor:more}(a), and part (d) of that corollary is the ratio $e_n/e_1$ for odd $n$. For $n=6$ the invariance of $e_2$, $e_4$, $e_6$, $e_3/e_1$ and $e_5/e_1$ is in \cite{KM}. By \eqref{eq:poly} the area determines the multiset of the tangent lengths. In fact it determines much more.

\begin{corollary}\label{cor:area}
Two polygons of the family with the same area either coincide or are symmetric to each other in the line $OI$. The family contains exactly two polygons that are symmetric in $OI$. The area is largest at one of them and smallest at the other, and it is strictly monotone on each of the two arcs of the family between them.
\end{corollary}

\begin{proof}
A polygon of the family, with its $n$ cyclic labellings, is a point of $\mathcal E_1$. As we saw, $S$ identifies the quotient of $\mathcal E_1$ by the reversal with a projective line. So two points of $\mathcal E_1$ with the same value of $S$ either coincide or are exchanged by the reversal. On real polygons the reversal is the reflection $\sigma$ in $OI$: the reflected polygon, labelled counterclockwise, has the tangent lengths $t_0,t_{-1},t_{-2},\dots$, and by Lemma \ref{lem:tt} two consecutive tangent lengths determine the polygon. This proves the first statement. The polygons of the family, without labels, form a circle, on which $\sigma$ acts as an involution that reverses the orientation. Such an involution has exactly two fixed points, and they divide the circle into two arcs exchanged by $\sigma$. By the first statement the continuous function $S$ is injective on each of the closed arcs, hence strictly monotone.
\end{proof}

For small $n$ the polygons of extreme area were known: see \cite{Radic04} for triangles, \cite{Josefsson,CM} for quadrilaterals and \cite{Radic05} for hexagons and octagons. Theorem \ref{thm:extremal} in Section \ref{ss:dist} tells which polygons these are for every $n$.

\begin{proposition}\label{prop:ten}
Over the Poncelet family
\begin{equation}\label{eq:ten}
\sum_{i=1}^n\bigl(|OB_i|^2-d^2\bigr)^5-10\,(R^2-d^2)\,r^4S^2
\end{equation}
is constant, and the coefficient $\beta$ in the formula for $\sum|OB_i|^8$ is
\begin{equation}\label{eq:beta}
\beta=-12d^4-8r^2\Bigl(R^2+d^2+e_2-\frac{e_3}{e_1}\Bigr).
\end{equation}
\end{proposition}

\begin{proof}
Recall that $|OB_i|^2-d^2=k-p_i$, and write $x=t_0$. Near $Q_+$ only $k-p_{-1}$ and $k-p_0$ have poles. Let $a$ be their half-sum and $b$ their half-difference. Since $t_{-1}$ and $t_1$ are the two roots of \eqref{eq:curve},
\[
t_{-1}+t_1=\frac{2kx}{x^2+r^2},\ \ \ t_{-1}t_1=\frac{r^2x^2+q^2}{x^2+r^2}.
\]
Hence $a=k-kx^2/(x^2+r^2)=kr^2/(x^2+r^2)$, while $b=x(t_1-t_{-1})/2$ behaves like $irx$. In
\[
(a+b)^5+(a-b)^5=2a^5+20a^3b^2+10ab^4
\]
only the last term has a pole, and it behaves like $10kr^6x^2$. Since $S^2$ behaves like $r^2x^2$, the sum of the fifth powers is an affine function of $S^2$ with the coefficient $10kr^4$. In the same way $(a+b)^3+(a-b)^3=2a^3+6ab^2$ has no pole, in agreement with the constancy of $\sum(|OB_i|^2-d^2)^3$.

For \eqref{eq:beta} we need one more term of the expansions. The formulas above give $p_{-1}+p_0=2k+O(x^{-2})$ and $p_{-1}p_0=r^2x^2+q^2-r^4+O(x^{-2})$. So the numbers $|OB_{-1}|^2=R^2-p_{-1}$ and $|OB_0|^2=R^2-p_0$ have the sum $2d^2+O(x^{-2})$ and the product $r^2x^2+R^4-2kR^2+q^2-r^4+O(x^{-2})$, and therefore
\[
|OB_{-1}|^8+|OB_0|^8=2r^4x^4+4r^2\bigl(R^4-2kR^2+q^2-r^4-4d^4\bigr)x^2+O(1).
\]
On the other hand $e_2=x\sum_{j\not\equiv0}t_j+e_2(t_j:j\not\equiv0)$ is constant, and at $Q_+$ the last term is the coefficient of $z^{n-3}$ in $\mathcal B(z)$, that is, $e_3/e_1$. Hence $S/r=x+(e_2-e_3/e_1)/x+O(x^{-2})$ and $S^4/r^4=x^4+4(e_2-e_3/e_1)x^2+O(x)$. Comparing the coefficients of $x^2$ in $\sum|OB_i|^8=2S^4+\beta S^2+\text{const}$ and using $q^2-r^4=k^2-2(R^2+d^2)r^2$ we get \eqref{eq:beta}.
\end{proof}

The curve $\mathcal E$ also shows that some quantities are never constant. Let us prove that $S^2+S^{*2}$ is not constant; by Theorem \ref{thm:limit}(d) the same is then true for $\sum\cos^2\varphi_i$. For even $n$ this is clear, since $S^*=S$ and $S$ has a pole at $Q_+$. Let $n$ be odd. The map $\kappa\colon(x,y)\mapsto(q/x,q/y)$ preserves \eqref{eq:curve} and commutes with $T$, and $S^*=S\circ\kappa$. An involution of $\mathcal E$ with a fixed point conjugates every translation to its inverse, and $T$ has order $n\ge3$. So $\kappa$ and $\iota$, which commute with $T$, are translations of order two, and so is $\kappa\iota$. If $S^*$ had a pole at $Q_+$, then $\kappa Q_+$ would be $T^jQ_+$ or $T^jQ_-=T^j\iota Q_+$ for some $j$, so that $\kappa=T^j$ or $\kappa\iota=T^j$. This is impossible for odd $n$, because $T^j$ has odd order. Hence $S^*$ is finite at $Q_+$, while $S$ has a pole there, and $S^2+S^{*2}$ is not constant.

\subsection{Products of distances}\label{ss:dist}
By Theorem \ref{thm:product} the product $\prod|XA_i|$ is constant for $X=I_m$. Here we find all points with this property. Let $d_m=|OI_m|$, so that $d_1=d$ and $\lambda_m^2=d_m/d$, and let $I_m'$ be the inverse of $I_m$ in $\G$; thus $|OI_m'|=d_m'=R^2/d_m$. Recall that $L$ and $L'$ are inverse in $\G$ as well.

\begin{theorem}\label{thm:dist}
\begin{enumerate}[label=\textup{(\alph*)}]
\item Let $X$ be the point of the line $OI$ with $\overrightarrow{OX}=(\xi/d)\,\overrightarrow{OI}$. Then
\begin{equation}\label{eq:dist}
\prod_{i=1}^n|XA_i|^2=K_X+\frac{\xi}{d\,r^2}\,\bigl((\xi-|OL|)(\xi-|OL'|)\bigr)^{\varepsilon}\prod_{1\le m<n/2}(\xi-d_m)^2(\xi-d_m')^2\cdot S^2,
\end{equation}
where $K_X$ does not depend on the polygon.
\item The product $\prod|XA_i|$ is constant over the Poncelet family if and only if $X$ is one of the points $O$, $I_m$, $I_m'$ with $1\le m<n/2$, or $n$ is even and $X$ is $L$ or $L'$.
\end{enumerate}
\end{theorem}

\begin{proof}[Proof of (a)]
For $\xi=0$ both sides are $R^{2n}$. Let $\xi\ne0$. Take the origin at $O$ and the first axis along $OI$. If $A$ is a point of $\G$ with the first coordinate $x$ and $t$ is the tangent length from $A$ to $\g$, then $|XA|^2=R^2+\xi^2-2\xi x$ and $t^2=R^2+d^2-r^2-2dx$. Hence
\[
|XA_j|^2=\frac{\xi}{d}\,(t_j^2-z^2),\ \ \ \text{where}\ \ z^2=R^2+d^2-r^2-d\Bigl(\xi+\frac{R^2}{\xi}\Bigr).
\]
By \eqref{eq:poly} with $s=S/r$,
\[
\prod_j(t_j^2-z^2)=\bigl(\mathcal A(z)+s\,\mathcal B(z)\bigr)\bigl(\mathcal A(-z)+s\,\mathcal B(-z)\bigr)=(-1)^n\bigl(\mathcal A(z)^2-s^2\mathcal B(z)^2\bigr),
\]
because $\mathcal A$ contains the powers $z^{n-j}$ with even $j$ and $\mathcal B$ those with odd $j$. By \eqref{eq:B},
\[
(-1)^{n+1}\mathcal B(z)^2=(-z^2)^{\varepsilon}\prod_{1\le m<n/2}\Bigl(-z^2-\frac{r_m^2}{\lambda_m^2}\Bigr)^2 .
\]
There are $n-1$ factors on the right. Multiply each of them by $\xi/d$ and use
\[
\frac{\xi}{d}\Bigl(-z^2-\frac{r_m^2}{\lambda_m^2}\Bigr)=(\xi-d_m)(\xi-d_m'),\ \ \ \frac{\xi}{d}\,(-z^2)=(\xi-|OL|)(\xi-|OL'|),
\]
which hold because $d_m+d_m'=2x_0+r_m^2/d_m$ by \eqref{eq:pencil} and $|OL|+|OL'|=2x_0$. This gives \eqref{eq:dist}.
\end{proof}

For the points of the line $OI$ part (b) follows from part (a): $S$ is not constant, and the coefficient of $S^2$ in \eqref{eq:dist} vanishes exactly at the points listed in (b). For the other points we need the vertices themselves and not only the tangent lengths. Let $R=1$, and identify the plane with $\mathbb C$ so that $O=0$ and $I=d$. Let $a_j$ be the complex number of the vertex $A_j$, and for a polynomial $H$ of degree at most $n$ put $H^*(z)=z^n\,\overline{H(1/\bar z)}$.

\begin{proposition}\label{prop:blaschke}
Let
\[
M(z)=z\,(z-|OL|)^{\varepsilon}\prod_{1\le m<n/2}(z-d_m)^2,\ \ \ \text{so that}\ \ \ M^*(z)=(1-|OL|\,z)^{\varepsilon}\prod_{1\le m<n/2}(1-d_mz)^2 .
\]
Then for every polygon of the family
\begin{equation}\label{eq:blaschke}
\prod_{j=1}^n(z-a_j)=M(z)+\eta\,M^*(z),\ \ \ \eta=(-1)^na_1a_2\cdots a_n .
\end{equation}
\end{proposition}

So the vertices are the solutions of the equation $M(z)/M^*(z)=-\eta$, and $M/M^*$ is a finite Blaschke product whose zeros are $O$, the centers $I_m$ taken twice, and $L$ if $n$ is even. The connection between Poncelet polygons inscribed in the unit circle and level sets of finite Blaschke products is well known \cite{DGSV}; for the case of Poncelet circles see \cite{ACG}. Here the zeros come out of the curve $\mathcal E$.

\begin{proof}
In the frame of Lemma \ref{lem:tt} the endpoints of the side with the unit normal $w=\cos\theta_0+i\sin\theta_0$ are $(r-it_0)w$ and $(r+it_1)w$, and that frame goes over into the present one by the half-turn $z\mapsto d-z$. So, by \eqref{eq:tt},
\[
a_0=d-(r-it_0)\,w,\ \ \ a_1=d-(r+it_1)\,w,\ \ \ w=\frac{t_0t_1-c+ir\,(t_0-t_1)}{2dr}.
\]
These formulas define $a_0$, and then $a_j=a_0\circ T^j$, as meromorphic functions on $\mathcal E$. For real polygons $|a_0|=1$ and $t_0^2=|a_0-d|^2-r^2$, hence
\[
t_0^2=1+d^2-r^2-d\,(a_0+a_0^{-1})
\]
on $\mathcal E$, and $a_0$ can have zeros and poles only at $Q_\pm$. At $Q_+$ we have $t_1=ir$, so $w$ behaves like $it_0/d$ and $a_0$ like $-t_0^2/d$: it has a double pole. Changing the signs of $t_0$ and $t_1$ in the formulas gives $\bar w$ and $\bar a_0$, and $\bar a_0=a_0^{-1}$ for real polygons. So $a_0\circ\iota=a_0^{-1}$ on $\mathcal E$, and $a_0$ has a double zero at $Q_-=\iota Q_+$. It follows that $a_j$ is finite at $Q_+$ for $j\not\equiv0$. Let $\alpha_j$ be its value there. We claim that
\begin{equation}\label{eq:alpha}
\alpha_m=\alpha_{n-m}=\frac{1}{d_m}\ \ (1\le m<n/2),\ \ \ \alpha_{n/2}=\frac{1}{|OL|}\ \ (n\text{ even}).
\end{equation}

Let $1\le m<n/2$. The chord $A_0A_m$ touches $\g_m$, and by Section \ref{ss:diag} the circle $\g_m$ lies on the left of the chord directed from $A_0$ to $A_m$. The formulas above hold for the pair $\G$, $\g_m$ as well, with $d_m$, $r_m$, $c_m=k_m-r_m^2$ and the tangent lengths $u_0=\lambda_mt_0$, $u_m=\lambda_mt_m$:
\[
a_0=d_m-(r_m-iu_0)\,w_m,\ \ \ a_m=d_m-(r_m+iu_m)\,w_m,\ \ \ w_m=\frac{u_0u_m-c_m+ir_m(u_0-u_m)}{2d_mr_m}.
\]
At $Q_+$ the function $u_m$ tends to $ir_m$ or to $-ir_m$, by \eqref{eq:taum}, and the relation \eqref{eq:curve} for $\g_m$ gives $u_m=\pm ir_m+k_m/u_0+\dots$ In the case of the minus sign the product $u_0(u_m+ir_m)$ would stay bounded, and so would $w_m$, and $a_0$ would have at most a simple pole. So $u_m=ir_m+k_m/u_0+\dots$ Thus $w_m$ behaves like $iu_0/d_m$ and $r_m+iu_m$ like $ik_m/u_0$, so that $\alpha_m=d_m+k_m/d_m=1/d_m$.

For the index $n-m$ we use the reflection $\sigma$ in $OI$. On real polygons it reverses the labels and replaces $a_j$ by $\bar a_{-j}=a_{-j}^{-1}$. So the reversal of the indices, as a map of $\mathcal E$, replaces $a_j$ by $a_{-j}^{-1}$, and its composition with $\iota$ replaces $a_j$ by $a_{-j}$. This composition fixes $Q_+$, because $t_{-1}=-ir$ there. Hence $\alpha_{-m}=\alpha_m$. Finally, for even $n$ the points $A_0$, $L$ and $A_{n/2}$ are collinear by Proposition \ref{prop:conj}, which for points of the unit circle and the real point $L$ means $a_{n/2}=(|OL|-a_0)/(1-|OL|\,a_0)$. At $Q_+$ this gives $\alpha_{n/2}=1/|OL|$, and \eqref{eq:alpha} is proved.

Now let $\mathcal V(z)=\prod_j(z-a_j)$, $\mathcal N(z)=\prod_{j\not\equiv0}(z-\alpha_j)$ and $g=\sum_ja_j$. The coefficients of $\mathcal V$ are invariant under the shifts, and their only poles are at the points of the orbit of $Q_+$, where they have order at most two. Near $Q_+$ we have
\[
\mathcal V(z)+g\,\mathcal N(z)=z\prod_{j\not\equiv0}(z-a_j)-a_0\Bigl(\prod_{j\not\equiv0}(z-a_j)-\mathcal N(z)\Bigr)+\mathcal N(z)\sum_{j\not\equiv0}a_j ,
\]
and the difference in the big brackets vanishes at $Q_+$. So the coefficients of $\mathcal V+g\,\mathcal N$ are functions on $\mathcal E_1$ with at most one simple pole. The curve $\mathcal E_1$ has genus one, being the quotient of $\mathcal E$ by a group of translations, so these functions are constant. Evaluating at $z=0$, where $\mathcal N(0)\ne0$, we can express $g$ through $\eta=\mathcal V(0)$, and we get
\[
\mathcal V(z)=\mathcal F(z)+\eta\,\frac{\mathcal N(z)}{\mathcal N(0)}=\mathcal F(z)+\eta\,M^*(z)
\]
with a polynomial $\mathcal F$ that does not depend on the polygon; the second equality is \eqref{eq:alpha}. For real polygons $|a_j|=1$, hence $|\eta|=1$ and $\mathcal V=\eta\,\mathcal V^*$, that is,
\[
\mathcal F+\eta\,M^*=M+\eta\,\mathcal F^* .
\]
The function $\eta$ has poles on $\mathcal E$, so it takes infinitely many values on the family. Therefore $\mathcal F=M$.
\end{proof}

\begin{proof}[End of the proof of Theorem \ref{thm:dist}]
Let $X$ be the point with the complex number $z$. By \eqref{eq:blaschke}, $\prod|XA_j|=|M(z)+\eta M^*(z)|$, and
\[
|M(z)+\eta M^*(z)|^2=|M(z)|^2+|M^*(z)|^2+2\operatorname{Re}\bigl(\eta\,M^*(z)\overline{M(z)}\bigr).
\]
This is constant if and only if $M(z)M^*(z)=0$, because $\eta$ takes infinitely many values on the unit circle. The zeros of $M$ are $O$, the points $I_m$ and, for even $n$, the point $L$; the zeros of $M^*$ are the inverse points $I_m'$ and, for even $n$, the point $L'$.
\end{proof}

Formula \eqref{eq:dist} and Proposition \ref{prop:blaschke} also give the polygons of extreme area. Since $S=rP/2$, where $P$ is the perimeter, they are the polygons of extreme perimeter as well.

\begin{theorem}\label{thm:extremal}
Let $A_-$ and $A_+$ be the points of $\G$ farthest from $I$ and nearest to $I$, and let $B_-$ be the point of $\g$ nearest to $O$.
\begin{enumerate}[label=\textup{(\alph*)}]
\item The area and the perimeter are largest for the polygon of the family with a vertex at $A_-$. Both tangent segments from this vertex have the largest possible length $\sqrt{(R+d)^2-r^2}$.
\item The area and the perimeter are smallest for the polygon of the family with a side touching $\g$ at $B_-$. For odd $n$ this polygon has a vertex at $A_+$; for even $n$ it has two sides perpendicular to $OI$.
\item The extreme values of the area satisfy
\begin{equation}\label{eq:gap}
S_{\max}^2-S_{\min}^2=4Rd\,r^2\Bigl(\frac{|OL|}{R}\Bigr)^{\varepsilon}\prod_{1\le m<n/2}\Bigl(\frac{|OI_m|}{R}\Bigr)^2 .
\end{equation}
\end{enumerate}
\end{theorem}

Put $\g$ with its center at $(0,0)$ and $\G$ with its center at $(a,0)$, $a>0$. Then $A_-=(R+a,0)$ and $B_-=(r,0)$, and (a) and (b) are the two conjectures at the end of \cite{CM}, where they are proved for the perimeter of quadrilaterals. For triangles see \cite{Radic04}, for the area of quadrilaterals \cite{Josefsson,Josefsson15}.

\begin{proof}
(a) Let $X=A_-$ in \eqref{eq:dist}, that is, $\xi=-R$. Then $\xi/(dr^2)<0$, the squared factors are positive, and for even $n$ so is $(\xi-|OL|)(\xi-|OL'|)$. So the coefficient of $S^2$ is negative, while the left-hand side is non-negative and vanishes exactly when $A_-$ is a vertex. Hence $S$ is largest exactly for the polygon with a vertex at $A_-$. The tangent length from a point of $\G$ is largest when the point is farthest from $I$, and $|IA_-|=R+d$.

(b) By Corollary \ref{cor:area} the smallest area belongs to the second polygon of the family that is symmetric in $OI$. If $\mathcal P$ is symmetric, then $\sigma(A_j)=A_{b-j}$ for some $b$, indices modulo $n$. A vertex fixed by $\sigma$ is $A_-$ or $A_+$, and a side fixed by $\sigma$ is perpendicular to $OI$ and touches $\g$ at $B_-$ or at the other point $B_+$ of $\g$ on $OI$. For odd $n$ there are exactly one fixed vertex and one fixed side. The polygon with a vertex at $A_+$ is symmetric, because its reflection is again a polygon of the family with a vertex at $A_+$, and it is not the polygon of (a), since it has only one vertex on $OI$. So it has the smallest area. The line $OI$ meets its boundary at $A_+$ and at the point of contact of the fixed side. As $I$ lies inside the polygon, these two points lie on different sides of $I$, so the point of contact is $B_-$. For even $n$ a symmetric polygon has either two fixed vertices and no fixed side, or two fixed sides and no fixed vertex, according to the parity of $b$. The polygon of (a) has the fixed vertex $A_-$, so the other symmetric polygon has two fixed sides, and they touch $\g$ at $B_-$ and $B_+$. In both cases a polygon with a side touching $\g$ at $B_-$ is unique, because this side determines the polygon.

(c) Let $R=1$ and let $C=(-|OL|)^{\varepsilon}\prod_md_m^2$ be the coefficient of $z^{n-1}$ in $M^*(z)$. Comparing the coefficients of $z^{n-1}$ in \eqref{eq:blaschke} we get $\sum_ja_j=b-\eta\,C$ with a real constant $b$. Since $t_j^2=1+d^2-r^2-2d\operatorname{Re}a_j$ and $e_2$ is constant,
\[
\frac{S^2}{r^2}=\sum_jt_j^2+2e_2=\text{const}+2dC\operatorname{Re}\eta .
\]
Putting $z=a_0$ in \eqref{eq:blaschke} gives $\eta=-M(a_0)/M^*(a_0)$. Here $M/M^*$ is a finite Blaschke product, and $a_0$ runs over the unit circle, so $\eta$ takes every value of modulus one. Hence $S_{\max}^2-S_{\min}^2=4d\,r^2|C|$, and \eqref{eq:gap} follows by scaling.
\end{proof}

\subsection{Tangents of the half-angles}
The statements of this subsection were first observed numerically. Let $T_i=\tan(\varphi_i/2)$, where $\varphi_i=\angle A_iLA_{i+1}$. One has the formula
\begin{equation}\label{eq:tanhalf}
\tan\frac{\varphi_i}{2}=\frac{2\,\dist(L,A_iA_{i+1})}{(\mu^2-1)\,|A_iA_{i+1}|},
\end{equation}
which follows from Lemma \ref{lem:casey} and the area of the triangle $A_iLA_{i+1}$.

\begin{theorem}\label{thm:tan}
For every $n$,
\begin{equation}\label{eq:tan1}
\tan\frac{\varphi_i}{2}=\frac{t_it_{i+1}+q}{r\,(t_i+t_{i+1})} .
\end{equation}
Let $n\ne4$. Then $q\ne r^2$,
\begin{equation}\label{eq:tan2}
\tan\psi_i=\frac{r}{r^2-q}\Bigl(t_{i+1}+\frac{q}{t_{i+1}}\Bigr),
\end{equation}
and the three sums
\[
\sum_{i=1}^nT_iT_{i+1},\ \ \ \sum_{i=1}^n\frac{1}{T_iT_{i+1}},\ \ \ \sum_{i=1}^n\Bigl(\frac{T_i}{T_{i+1}}+\frac{T_{i+1}}{T_i}\Bigr)
\]
are constant over the Poncelet family. The same is true for the three sums obtained by replacing each $T_i$ with $\tan\psi_i$.
\end{theorem}

The case $n=4$ is different. There $\varphi_i+\varphi_{i+1}=\pi$, so $T_{i+1}=1/T_i$; the first two sums are equal to $4$, the third is $4(T_1^2+T_1^{-2})$ and is not constant, and $\tan\psi_i$ is not defined. Indeed, for the polygon with the vertex $A_1$ on $OI$ the symmetry gives $\varphi_1=\varphi_4$, hence $\varphi_1=\pi/2$ and $T_1=1$. The two sides perpendicular to $OI$ of the other symmetric polygon are chords at the distances $|d-r|$ and $d+r$ from $O$, with the tangent lengths $t$ and $q/t=r^2/t$, so that $t\ne r$ for one of them. If $A_1A_2$ is this side, then $t_1=t_2=t$ and \eqref{eq:tanhalf} gives $T_1=(t^2+r^2)/(2rt)\ne1$.

\begin{proof}
\emph{The formulas.} In the frame of Lemma \ref{lem:tt} the point $L$ is $(-|IL|,0)$, so $\dist(L,A_iA_{i+1})=r+|IL|\cos\theta_i$, and $2dr\cos\theta_i=p_i-c$. The points $L$ and $L'$ are inverse in $\g$, so $|IL|\cdot|IL'|=r^2$, and $2d\,|IL'|-c=2d(|OL'|-x_0)=q$. Hence
\[
\dist(L,A_iA_{i+1})=\frac{|IL|}{2dr}\,(p_i+q),
\]
and \eqref{eq:tanhalf} with $\mu^2-1=|IL|/d$ gives \eqref{eq:tan1}.

Fix $i$ and write $x=t_{i+1}$, $y=t_i$, $z=t_{i+2}$. By Lemma \ref{lem:vieta}, $y$ and $z$ are the two roots of \eqref{eq:curve}, so
\begin{equation}\label{eq:yz}
y+z=\frac{2kx}{x^2+r^2},\ \ \ yz=\frac{r^2x^2+q^2}{x^2+r^2}.
\end{equation}
It follows that $(x+y)(x+z)=x^2Z/(x^2+r^2)$, where $Z=x^2+q^2/x^2+2(k+r^2)$. Inserting \eqref{eq:yz} into \eqref{eq:tan1} we find
\begin{equation}\label{eq:TT}
T_i-T_{i+1}=\frac{(y-z)(x^2-q)}{r\,(x+y)(x+z)},\ \ \ T_i+T_{i+1}=\frac{2(x^2+q)(k+r^2+q)}{r\,x\,Z},\ \ \ T_iT_{i+1}=1+\frac{C}{Z},
\end{equation}
with $C=2(q-r^2)(k+r^2+q)/r^2$. If $q=r^2$, then $T_iT_{i+1}=1$, that is, $\varphi_i+\varphi_{i+1}=\pi$ for all $i$, and $\sum\varphi_i=2\pi$ gives $n=4$. So $q\ne r^2$ for $n\ne4$. Since $\psi_i=(\varphi_i+\varphi_{i+1})/2$, we have $\tan\psi_i=(T_i+T_{i+1})/(1-T_iT_{i+1})$, and \eqref{eq:TT} gives \eqref{eq:tan2}.

\emph{The mechanism.} Let $d\tau=d\theta_i/(t_i+t_{i+1})$ be the common differential of \eqref{eq:poncelet-diff}. Suppose that non-zero numbers $V_i$ attached to the polygon satisfy
\begin{equation}\label{eq:mech}
d(V_iV_{i+1})=a\,(V_i^2-V_{i+1}^2)\,d\tau
\end{equation}
with a constant $a$ that does not depend on $i$. Then $\sum V_iV_{i+1}$ is constant, and so is $\sum1/(V_iV_{i+1})$, because \eqref{eq:mech} implies $d\bigl(1/(V_iV_{i+1})\bigr)=a\,(V_i^{-2}-V_{i+1}^{-2})\,d\tau$. The first example is $V_i=t_i$: by \eqref{eq:tt},
\begin{equation}\label{eq:dp}
dp_i=-2dr\sin\theta_i\,d\theta_i=-r\,(t_i^2-t_{i+1}^2)\,d\tau .
\end{equation}
Differentiating also the first relation of \eqref{eq:tt} and solving for $dt_{i+1}$, we get $r\,dt_{i+1}=\bigl(kt_{i+1}-t_i(t_{i+1}^2+r^2)\bigr)d\tau$, which by \eqref{eq:yz} is
\begin{equation}\label{eq:dx}
dx=\frac{(x^2+r^2)(z-y)}{2r}\,d\tau .
\end{equation}

\emph{The sums for $T_i$.} By \eqref{eq:TT} and \eqref{eq:dx},
\[
d(T_iT_{i+1})=-\frac{C}{Z^2}\Bigl(2x-\frac{2q^2}{x^3}\Bigr)dx=\frac{C\,(x^4-q^2)(x^2+r^2)(y-z)}{r\,x^3Z^2}\,d\tau ,
\]
and, again by \eqref{eq:TT},
\[
T_i^2-T_{i+1}^2=\frac{2(y-z)(x^4-q^2)(k+r^2+q)(x^2+r^2)}{r^2x^3Z^2} .
\]
Comparing, $d(T_iT_{i+1})=\frac{q-r^2}{r}\,(T_i^2-T_{i+1}^2)\,d\tau$. So \eqref{eq:mech} holds for $V_i=T_i$, and the first two sums are constant. For the third put $P_i=T_iT_{i+1}$. By \eqref{eq:tan2} and \eqref{eq:TT}, $(x+q/x)^2=Z-2(k+r^2-q)$ and $Z=C/(P_i-1)$, so
\[
(T_i+T_{i+1})^2=(1-P_i)^2\tan^2\psi_i=\frac{r^2}{(r^2-q)^2}\Bigl(C\,(P_i-1)-2(k+r^2-q)(P_i-1)^2\Bigr).
\]
Dividing by $P_i$ we see that $T_i/T_{i+1}+T_{i+1}/T_i$ is a combination of $1$, $P_i$ and $1/P_i$ with constant coefficients, and the third sum is constant.

\emph{The sums for $\tan\psi_i$.} By \eqref{eq:tan2} it is enough to consider $v_i=t_i+q/t_i$. Put $z_i=p_i+q^2/p_i$. By \eqref{eq:p2},
\[
v_iv_{i+1}=z_i+q\Bigl(\frac{t_i}{t_{i+1}}+\frac{t_{i+1}}{t_i}\Bigr)=\Bigl(1-\frac{q}{r^2}\Bigr)z_i+\frac{2kq}{r^2},\ \ \
v_i^2+v_{i+1}^2=\frac{z_i\,(2k-z_i)}{r^2}+4q .
\]
Since $v_i^2-v_{i+1}^2=(t_i^2-t_{i+1}^2)(1-q^2/p_i^2)$, formula \eqref{eq:dp} gives
\[
d(v_iv_{i+1})=\Bigl(1-\frac{q}{r^2}\Bigr)\Bigl(1-\frac{q^2}{p_i^2}\Bigr)dp_i=-r\Bigl(1-\frac{q}{r^2}\Bigr)(v_i^2-v_{i+1}^2)\,d\tau .
\]
So \eqref{eq:mech} holds for $V_i=v_i$. Moreover, since $q\ne r^2$, the second relation above expresses $v_i^2+v_{i+1}^2$ as a quadratic polynomial in $v_iv_{i+1}$; dividing by $v_iv_{i+1}$ we obtain the third sum.
\end{proof}

Formula \eqref{eq:tan2} can be written as $\cos^2\psi_i=(r^2-q)^2/\bigl(|IA_{i+1}|^2\,|IA_{i+1}^*|^2\bigr)$, where $A_{i+1}^*$ is the vertex of the conjugate polygon; for $n=4$ both sides vanish. Since $\prod|IA_i^*|=\prod|IA_i|$, this gives the product formula of Theorem \ref{thm:limit}(e) once more, and shows that its right-hand side is $|r^2-q|^n$. At $L'$ the same computation, with $\mu'$ in place of $\mu$, gives
\[
\cos^2\psi_i'=\frac{(r^2+q)^2}{|IA_{i+1}|^2\,|IA_{i+1}^*|^2},
\]
because $k^2-4R^2r^2/\mu'^2=(r^2+q)^2$. So the constant in the product formula at $L'$ is $(r^2+q)^n$, also for $n=4$, although for $n=4$ the sum $\sum\cos^2\psi_i'$ is not constant.

In billiard terms $\varphi_i/2$ is the angle between the trajectory and the ellipse at the point of impact $P_i$; these angles are studied in \cite{Stachel}. So Theorem \ref{thm:tan} says that for periodic billiard trajectories of winding number one in an ellipse, with $n\ne4$, the sum of the products of the tangents of these angles at consecutive impacts is constant, and so is the sum of the products of their cotangents. There is another way to look at the numbers $T_i$, which comes from a construction used in \cite{AST} for billiards. The unit vectors in the directions $LA_i$ form a Poncelet polygon inscribed in the unit circle centered at $L$ and circumscribed about a concentric ellipse, and the $T_i$ are the tangent lengths of its polar polygon, which is circumscribed about the unit circle. So the $T_i$ play, for another Poncelet family, the role that the $t_i$ play for ours.

\begin{remark}\label{rem:halfprod}
The law of cosines in the triangle $A_iLA_{i+1}$ gives
\[
\cos^2\frac{\varphi_i}{2}=\frac{(\mu^2-1)(t_i+t_{i+1})^2}{4\mu^2\,t_it_{i+1}} .
\]
So, by \eqref{eq:odd-product}, the product $\prod\cos(\varphi_i/2)$ is constant for odd $n$. For even $n$ it is never constant. Indeed, then $\prod t_i=q^{n/2}$, and we claim that $\prod(t_i+t_{i+1})$ has a pole of order two at the point $Q_+$ of Section \ref{ss:curve}. The factors $t_{-1}+t_0$ and $t_0+t_1$ have simple poles there, and the other factors are finite. Suppose that $t_j+t_{j+1}$ vanishes at $Q_+$. Then the point $T^jQ_+$ has the form $(x,-x)$, so it is fixed by the map $\vartheta\colon(x,y)\mapsto(-y,-x)$ of $\mathcal E$. This map satisfies $\vartheta T^j=T^{-j}\vartheta$ and $\vartheta Q_+=T^{-1}Q_+$, hence $T^jQ_+=T^{-j-1}Q_+$ and $n$ divides $2j+1$, which is impossible for even $n$.
\end{remark}

\subsection{The constants \texorpdfstring{$\lambda_m$}{lambda m}}
Poncelet's general theorem is true for open chains as well. Let $X_0$ be a point of $\G$ and $X_j=f^j(X_0)$, so that the chain always uses the same branch of tangents. If $f^m\ne\mathrm{id}$, then, as $X_0$ moves, the line $X_0X_m$ remains tangent to a conic of the pencil \cite[Theorem 4.1.1]{BKOR}. If also $f^{2m}\ne\mathrm{id}$, the argument of Section \ref{ss:diag}, which does not use closure, shows that this conic is a circle $\g_m$ inside $\G$ containing $L$. So the circles $\g_m$ and the numbers $\lambda_m$ belong to the pair of circles and to $m$, whether polygons close or not. The coincidence of two of them forces closure. For instance, if $\g_2$ and $\g_3$ are the same circle, then $f^2(X)$ and $f^3(X)$ are joined to $X$ by tangents from $X$ to this circle, and so is $f^{-2}(X)$. Since $f$ has no fixed points, $f^3(X)\ne f^2(X)$, so $f^3=f^{-2}$ and $f^5=\mathrm{id}$. The polygons need not be convex: if $\g$ carries convex pentagons, then for the pair $\G$, $\g_2$ both circles are $\g$, and the polygons are pentagrams.

For a pair of circles that carries $n$-gons the constants $\lambda_m$ can be computed from the curve $\mathcal E$. Let $\tau_j$ be the value of $t_j$ at $Q_+$, as in the proof of Theorem \ref{thm:poly}. Then $\tau_1=ir$, and \eqref{eq:curve} and Lemma \ref{lem:vieta} give $\tau_j=i\sigma_j$ with real numbers
\begin{equation}\label{eq:sigma}
\sigma_1=r,\ \ \ \sigma_2=\frac{r^4-q^2}{2kr},\ \ \ \sigma_{j+1}=\frac{2k\,\sigma_j}{r^2-\sigma_j^2}-\sigma_{j-1}\ \ (2\le j\le n-2).
\end{equation}
The denominators do not vanish: if $\tau_j=\pm ir$, then the coefficient of $y^2$ in \eqref{eq:curve} vanishes at $x=\tau_j$, so one of $\tau_{j-1}$, $\tau_{j+1}$ is infinite, while $t_{j\pm1}$ is finite at $Q_+$ for $1\le j\pm1\le n-1$.

Let $1\le m<n/2$. By \eqref{eq:taum}, $r_m^2/\lambda_m^2=-\tau_m^2=\sigma_m^2$. On the other hand $r_m^2=\lambda_m^4d^2-2\lambda_m^2dx_0+R^2$ by \eqref{eq:pencil}. Therefore
\begin{equation}\label{eq:lambda-m}
\lambda_m^2=\frac{\sigma_m^2+2dx_0-\sqrt{(\sigma_m^2+2dx_0)^2-4d^2R^2}}{2d^2},\ \ \ 2dx_0=R^2+d^2-r^2 ,
\end{equation}
where the smaller root is taken because $|OI_m|<|OL|$. For $m=2$ this gives again $\lambda_2=2Rr/k$. Moreover, the relation $t_{-j}=-t_j$ at $Q_+$ means that $\sigma_{n-j}=-\sigma_j$. In particular
\begin{equation}\label{eq:closure}
\sigma_{n/2}=0\ \ \text{for even }n,\ \ \ \ \sigma_{(n-1)/2}+\sigma_{(n+1)/2}=0\ \ \text{for odd }n .
\end{equation}
These are conditions on $R$, $r$, $d$ satisfied by every pair of circles that carries $n$-gons. For $n=3$ and $n=4$ they read $R^2-d^2=2Rr$ and $q=r^2$, which are the formulas of Euler and Fuss. The recursion also computes the invariants of Theorem \ref{thm:poly}, because $\mathcal B(z)=z^{\varepsilon}\prod_{1\le m<n/2}(z^2+\sigma_m^2)$.

The conditions \eqref{eq:closure} can also be used as a test of closure. Let $\G$, $\g$ be any two circles with $d+r<R$. The curve $\mathcal E$, the translation $T$ and the point $Q_+$ are defined as before, and \eqref{eq:sigma} computes $\tau_j=i\sigma_j$ as long as the values are finite. Conditions for closure that avoid elliptic functions were also given by Cie\'slak and Szczygielska \cite{CS10}.

\begin{proposition}\label{prop:test}
Suppose that $\sigma_1,\dots,\sigma_{\lceil n/2\rceil}$ are defined by \eqref{eq:sigma}, are finite and satisfy \eqref{eq:closure}. Then the Poncelet map $f$ has period exactly $n$. The polygons are convex if and only if $\sigma_j>0$ for $1\le j<n/2$.
\end{proposition}

\begin{proof}
Let $\mathcal R$ be the reversal of the indices, which sends the point $(t_0,t_1)$ of $\mathcal E$ to $(t_0,t_{-1})$, and let $\iota(x,y)=(-x,-y)$. The map $\iota\mathcal R$ fixes $Q_+$, because $t_{-1}=-ir$ at $Q_+$, and $(\iota\mathcal R)T^j=T^{-j}(\iota\mathcal R)$. Let $n=2N$. Then $\sigma_N=0$ means that $T^NQ_+$ has the first coordinate $0$. The two roots of \eqref{eq:curve} with $x=0$ are opposite, so $\iota\mathcal R$ fixes the points of $\mathcal E$ with $x=0$. Hence $T^NQ_+=\iota\mathcal R\,T^NQ_+=T^{-N}Q_+$. Let $n=2N+1$. Then $T^NQ_+=(x,-x)$ for some $x$, and with the map $\vartheta$ of Remark \ref{rem:halfprod} we get $T^NQ_+=\vartheta T^NQ_+=T^{-N}\vartheta Q_+=T^{-N-1}Q_+$. In both cases the translation $T^n$ has a fixed point, so $T^n=\mathrm{id}$. If the order $b$ of $T$ were a proper divisor of $n$, then $b\le N$, $T^bQ_+=Q_+$, and $t_b$ would have a pole at $Q_+$, contrary to the finiteness of $\sigma_b$. So $T$ has order $n$, and so has $f$, since a translation of $\mathcal E$ that fixes the real points is the identity.

Thus every polygon closes after $n$ steps and goes $w$ times around $\G$, where $\gcd(w,n)=1$, and $f$ moves each vertex $w$ places forward in the cyclic order of the vertices on $\G$. The chord directed from $X$ to $f(X)$ has $\g$, and so $L$, on its left. So the angle of the counterclockwise rotation about $L$ from $LX$ to $Lf(X)$ is less than $\pi$, and $2w\pi<n\pi$. Let $1\le j<n/2$. Then $f^j\ne\mathrm{id}$ and $f^{2j}\ne\mathrm{id}$, so the chords $Xf^j(X)$ touch a circle $\g_j$ inside $\G$ containing $L$. The argument of Section \ref{ss:diag} now gives $\sum\alpha(A_i)=2\pi\,(jw\bmod n)$, so $\g_j$ lies on the left of the chord directed from $A_0$ to $A_j$ if $jw\bmod n<n/2$ and on its right if $jw\bmod n>n/2$. In the first case the proof of Proposition \ref{prop:blaschke} gives $\lambda_jt_j(Q_+)=ir_j$, that is, $\sigma_j>0$. In the second case the same argument applied to the chord directed from $A_j$ to $A_0$ gives $\lambda_jt_j(Q_+)=-ir_j$, that is, $\sigma_j<0$. If $w=1$, then $jw<n/2$ for all $j<n/2$. If $w\ge2$, then $n\ge5$, and $j=\lfloor n/(2w)\rfloor+1$ satisfies $j<n/2$ and $n/2<jw<n$. So all $\sigma_j$ with $j<n/2$ are positive if and only if $w=1$, that is, if and only if the polygons are convex.
\end{proof}

For example, let $R=1$ and $d=0.35$. For $n=5$ the condition $\sigma_2+\sigma_3=0$ has two solutions, $r\approx0.631902$ with $\sigma_2\approx0.113999$ and $r\approx0.277186$ with $\sigma_2\approx-1.228293$. The first pair of circles carries convex pentagons, the second one pentagrams.

\begin{question}
Find closed formulas for the constants $\lambda_m$, $m\ge3$, in terms of $R$, $r$, $d$, and for the conditions \eqref{eq:closure}.
\end{question}

\section{Concluding remarks}

The tangent lengths $t_i$ are good coordinates on the Poncelet family between two circles. Three facts carry the whole paper: the relation \eqref{eq:Phi} between neighbours, Casey's Lemma \ref{lem:casey}, which makes the distances to the limiting points and the tangent lengths to all circles of the pencil proportional to the $t_i$, and the constancy of the sums $\sum_it_it_{i+m}$, $1\le m\le n-1$. The first two are identities, and the third is Weill's theorem applied to each circle $\g_m$ of Poncelet's general theorem.

We do not think that a long list of invariants is a goal in itself. There are infinitely many of them, as Proposition \ref{prop:liouville} shows, and most are combinations of a few. What we find worth recording is, first, that the few basic ones have proofs within school geometry, Section \ref{s:nocalc}; second, that the invariants of degree four and six involving the area have universal coefficients; and third, that the theorems on the sum and on the product of cosines for billiards in ellipses are consequences of the theorem that Weill proved in 1878.

It would be interesting to have elementary proofs of Corollary \ref{cor:more} and of Theorems \ref{thm:poly}, \ref{thm:dist} and \ref{thm:extremal}, and proofs of Theorems \ref{thm:product} and \ref{thm:tan} without derivatives.

\section*{Acknowledgments}
The author thanks Claude Fable 5.2 and GPT 6.1 Sol for their help with the English language and the presentation of the paper. The mathematical content is entirely the author's responsibility.

\end{document}